\documentclass[12pt,oneside]{article}
\usepackage{epsfig}
\usepackage{latexsym}
\usepackage{amsmath}
\usepackage{amssymb}
\usepackage[english]{babel}
\usepackage{xcolor}

\newtheorem{thm}{Theorem}[subsection]
\newtheorem{cor}[thm]{Corollary}
\newtheorem{lem}[thm]{Lemma}
\newtheorem{prop}[thm]{Proposition}
\newtheorem{defn}[thm]{Definition}

\newtheorem{expl}[thm]{Example}

\newcommand{\Real}{\mathbb{R}}

\numberwithin{equation}{section}

\newfont{\sfl}{cmssi12}

\begin{document}

\title{ Quasi-Leontief Utility Functions on Partially Ordered Sets\, I: Efficient Points.}

\author{\thanks {University of Perpignan, Department of Economics, 52 avenue Paul Alduy,
66800 Perpignan, France.}  Walter Briec,  QiBin Liang and Charles
Horvath \thanks {University of Perpignan, 52 avenue Paul Alduy,
66800 Perpignan, France.}}  \maketitle \
\\


\bigskip
\textbf{Abstract:} A function $u: X\to\mathbb{R}$ defined on a partially ordered set is quasi-Leontief if, if for all $x\in X$, the upper level set $\{x^\prime\in X: u(x^\prime)\geqslant u(x)\} $ has a smallest element. A function 
$u: \prod_{j=1}^nX_j\to\mathbb{R}$ whose partial functions obtained by freezing $n-1$ of the variables are all quasi-Leontief is an individually quasi-Leontief function; a point $x$ of the product space is an efficient point for $u$ if it is a minimal element of  $\{x^\prime\in X: u(x^\prime)\geqslant u(x)\} $. 
Part I deals with the maximisation of quasi-Leontief functions and the existence of efficient maximizers.    Part II is concerned with the existence of efficient Nash equilibria for abstract games whose payoff functions are individually quasi-Leontief.   Order theoretical and algebraic arguments are dominant in the first part while, in the second part, topology is heavily involved.   
In the framework and the language of tropical algebras, our quasi-Leontief functions are the additive functions defined on a semimodule with values in the semiring of scalars.  

\newpage
\section{Introduction}\label{secintro}
A function $u : X\to\Lambda$ defined on a partially ordered set $X$ with values in a totally ordered set $\Lambda$ is a quasi-Leontief function if for all $x\in X$ the upper level set $\{x^\prime\in X: u(x^\prime)\geqslant u(x)\} $ has a smallest element $x^\circ$;   if $x = x^\circ$ then $x$ is an efficient point for $u$.  If $X = \mathbb{R}^n_+$  and $\Lambda =  \mathbb{R}_+$ then the classical Leontief functions are exactly the positively homogeneous quasi-Leontief functions.

\medskip\noindent
This paper can  be interpreted as a contribution to monotone analysis in the framework of the so called tropical algebras or Maslov semilattices, \cite{agk} and the references therein. In this framework, writing as usual $x\oplus y$ for the underlying idempotent operation of the algebra in question, quasi-Leontief functions are those functions for which $u(x\oplus y) = u(x) \oplus u(y)$; they are therefore the additive maps of the tropical algebra in question. We have already  said that, for $X = \mathbb{R}_+^n$ and $\Lambda = \mathbb{R}_+$, the classical Leontief functions are exactly the homogeneous quasi-Leontief  functions;  $\mathbb{R}_+^n$ equipped with its inf-semilattice structure is a tropical algebra whose semiring of scalars is $\mathbb{R}_+$, 
multiplication by scalars is here the usual multiplication. Leontief functions are therefore those functions $u : \mathbb{R}_+^n\to \mathbb{R}_+$ such that $u(tx\oplus y) = tu(x) \oplus u(y)$: in the context of tropical algebras, the Leontief functions are the linear maps.

\medskip\noindent
Section \ref{secqf} deals with quasi-Leontief utility functions and their basic properties. To each quasi-Leontief function $u : X\to\Lambda$ is associated a dual map 
$u^\sharp : u(X)\to X$ such that $u(x)\geqslant t$ if and only if 
$x\geqslant u^\sharp(t)$. Much of the properties of quasi-Leontief functions follow from this duality. Quasi-leontief functions are characterized in Sections \ref{subsecarqf} and  \ref{subsecsemilat}.
If $X$ is a closed convex subset of   $ \mathbb{R}^n$ then a quasi-Leontief function on $X$ is quasi-convex and upper semicontinuous; but one has to notice that the domain $X$ is only assumed to be partially ordered, and can therefore be far from being convex even if it is a subset of a vector space. If $X$ is an inf-semillatice -- any two elements of $X$ have a greatest lower bound -- then quasi-Leontief functions are inf-preserving functions. 
Section \ref{subsecmaxql} deals with maximisation of quasi-Leontief functions on a given subset $S$ of $X$. \\
A function $u$ defined on a product of partially ordered sets $ X_1\times\cdots\times X_n$ is globally quasi-Leontief if it is quasi-Leontief with respect to the product partial order; it is individually quasi-Leontief if, $n-1$ of the variables being frozen, it is quasi-Leontief in the remaining variable. 
An efficient point $x^\flat$ for an individually quasi-Leontief function $u$ is a minimal element of $\{x^\prime\in X_1\times\cdots\times X_n: u(x^\prime)\geqslant u(x^\flat)\}$; we show that, if $\arg\!\max(u, S)\neq\emptyset$, where $S$ is a product of comprehensive subsets $S_i\subset X_i$, then for all $x^\star\in\arg\!\max(u, S)$ there is an efficient point  
$x^\flat\in\arg\!\max(u, S)$  such that $x^\star\geqslant x^\flat$, Theorem \ref{existseffglob}.

\vspace{1cm}
\section{Quasi-Leontief utility functions}\label{secqf}
All the sets under consideration are partially ordered sets; the relevant examples are $\mathbb{R}$ with its usual ordering, $\mathbb{R}^n$ with the partial order associated to the positive cone 
$\mathbb{R}^n_+ = \{(x_1, \cdots, x_n) : \forall i \, x_i\geqslant 0\}$, Riesz spaces, $\mathbb{R}^n$ equipped with the lexicographic order or subsets thereof equipped with the induced partial order. A single notation will be used for all the partial orders under consideration. 

\bigskip

 We recall that a function $u : L^\prime\to L$ from a partially ordered set $L^\prime$ to a partially ordered set $L$  is {\bf isotone} if $x^\prime\geqslant x$ implies $u(x^\prime)\geqslant u(x)$. Given a partially ordered set $(L, \leqslant)$ and an element $x\in L$ we set $\downarrow\!\!(x) = \{x^\prime\in L: x^\prime\leqslant x\}$ and $\uparrow\!\!(x) = \{x^\prime\in L: x^\prime\geqslant x\}$; given two elements $x_1$ and $x_2$ the interval $[x_1, x_2]$ is the set  $\uparrow\!\!(x_1)\, \cap\downarrow\!\!(x_2)$, which can be empty. Given  a subset $S$ of $L$ let ${\displaystyle\downarrow\!\!(S) = \cup_{x\in S}\downarrow\!\!(x)}$ and 
$\uparrow\!\!(S) = \cup_{x\in S}\uparrow\!\!(x)$;  A subset $S$ of $L$   is upward (respectively downward) if  
 $S = \uparrow\!\!(S)$ (respectively $S = \downarrow\!\!(S)$). Said differently, 
  $S = \downarrow\!\!(S)$ if $x^\prime\in S$ and $x^\prime\geqslant x$ implies 
  $x\in S$.
  Downward sets will be called {\bf comprehensive}. \\
  Given two comparable elements $x_1$ and $x_2$, let us say $x_1\geqslant x_2$, of a partially ordered set 
  $L$ the {\bf interval $\boldsymbol{[x_1, x_2]}$} is the set 
  $\{x\in L: x_1\leqslant x \leqslant x_2\}$.
 
 \medskip\noindent
 The partially ordered set $L$ is totally ordered if, for all $x_1$, $x_2$ in $L$ either $x_1\leqslant x_2$ or  $x_2\leqslant x_1$, it will be convenient to count the empty set among the totally ordered sets ; a nonempty subset $E$ of $L$ which is  totally ordered with respect to the induced partial order is an order {\bf chain} in $L$. The set of real numbers is totally ordered as well as ${\mathbb R}^n$ with the lexicographic order.  \\
An element $\bar{x}\in L$ is a least (respectively largest) element of a subset $S\subset L$ if, for all $x\in S$, $\bar{x}\leqslant x$ (respetively $\bar{x}\geqslant x$ ) and $\bar{x}\in S$; a given subset $S\subset L$ as at most one least (respectively largest)  element.

 \subsection{Efficient points and quasi-Lontief functions}\label{subsecefficient}
 
 \medskip\noindent . 
 The general framework, that will be kept constant throughout this section and most of the paper, is entirely algebraic : a function $u : X\to\Lambda$, ``the utility function'', where   $(X, \leqslant)$ is a partially ordered set and $(\Lambda, \leqslant)$ is  a totally ordered set; $X$ could be a subset of $\mathbb{R}^n$ partially ordered by the positive cone $\mathbb{R}^n_+$, or $\mathbb{Z}^n$  for example;    and $\Lambda$ could be a subset of 
 $\mathbb{R}$ or $\mathbb{R}^n$ with the lexicographic ordering. Even if an interval of $\mathbb{R}$ is the natural choice for $\Lambda$ none of the specificities of intervals of $\mathbb{R}$ comes into play for most of our results; so we keep the ``general $\Lambda$'' and we add specific requirements when and where they are needed.\\
 When we will come to Nash equilibria some topological considerations, connectedness for example, will exclude spaces such as  $\mathbb{Z}^n$ from our considerations.
 
\begin{defn}\label{defefficient}
Given a utility function $u : X\to\Lambda$ we say that  
$x\in X$ is efficient, with respect to $u$, if, for all $x^\prime\in X$, $u(x^\prime)\geqslant u(x)$ implies 
$x^\prime\geqslant x$. In other words, $x$ is efficient if $x = \min\{x^\prime\in X:  u(x^\prime)\geqslant u(x)\}$. Given a subset $S$ of $X$ we will denote by ${\mathcal E}(u; S)$ the set, possibly empty, of efficient points which belong to $S$.
\end{defn}

\noindent By  ${\mathcal E}(u; S)$ me mean, as the definition says,  
${\mathcal E}(u; X)\cap S$ and not the set of efficient points of the restriction of $u$ to $S$, that is ${\mathcal E}(u_{\vert S}; S)$ ; these sets do not have to be the same as one can see by taking $S$ to be a singleton. It is clear from the definition that an element of 
${\mathcal E}(u; S)$ is also an efficient point of the restriction $u_{\vert S}: S\to\Lambda$.

\begin{lem}\label{efficient1} Assume that  $u : X\to\Lambda$ is isotone. Then, for all $x^\circ\in X$ the following statements are equivalent: \\
$(1)$ $x^\circ\in X$ is efficient; \\
$(2)$ for all $x\in X$, $u(x)\geqslant u(x^\circ)$ if and only if $x\geqslant x^\circ$; \\
$(3)$ $\uparrow\!\!(x^\circ) = u^{-1}\big(\uparrow\!\!u(x^\circ)\big)$  .
\end{lem}
{\it Proof:} Given that $u$ is isotone the equivalence of $(1)$ and $(2)$ is clear and $(3)$ is exactly $(2)$.  

\begin{lem}\label{efficientchain} If  $u : X\to\Lambda$ is isotone then ${\mathcal E}(u; X)$ is a totally ordered subset of X.
\end{lem}
{\it Proof:} Let $x_1^\circ$ and $x_2^\circ$ be to efficient points. Since $\Lambda$ is totally ordered we have either 
$u(x_1^\circ)\geqslant u(x_2^\circ)$ or $u(x_2^\circ)\geqslant u(x_1^\circ)$ and therefore, by $(2)$ of Lemma \ref{efficient1}, 
either 
$x_1^\circ\geqslant x_2^\circ$ or $x_2^\circ\geqslant x_1^\circ$. \hfill$\Box$ 

\begin{lem}\label{efficient2} Assume that  $u : X\to\Lambda$ is isotone. Then, for all $x\in X$ there exists at most one efficient point 
$x^\circ$ such that  $u(x) = u(x^\circ)$ and, assuming  that efficient point exists, $x^\circ = \min\{x^\prime\in X: u(x^\prime)\geqslant u(x)\}$. Reciprocally, if  $x^\circ = \min\{x^\prime\in X: u(x^\prime)\geqslant u(x)\}$ then $x^\circ$ is efficient and $u(x) = u(x^\circ)$ and therefore,
$x^\circ = \min\{x^\prime\in X: u(x^\prime) = u(x)\}$ .
\end{lem}
{\it Proof:} If $u(x) = u(x^\circ_{i})$, $i = 1, 2$, with $x^\circ_{i}$ efficient,  then, from $u(x^\circ_{1})\geqslant u(x^\circ_{2})$ and 
$u(x^\circ_{2})\geqslant u(x^\circ_{1})$ and $(2)$ of Lemma \ref{efficient1} we have $x^\circ_{1} = x^\circ_{2}$.\\
If $u(x^\prime)\geqslant u(x)$ and $u(x) = u(x^\circ)$ then  $u(x^\prime)\geqslant u(x^\circ)$ and therefore $x^\prime\geqslant x^\circ$; this proves that $x^\circ = \min\{x^\prime\in X: u(x^\prime)\geqslant u(x)\}$.

\medskip\noindent
Assume that $x^\circ = \min\{x^\prime\in X: u(x^\prime)\geqslant u(x)\}$.  From the trivial inequality 
$u(x)\geqslant u(x)$ and the minimality of $x^\circ$ we have $x\geqslant x^\circ$, and since $u$ is isotone, 
$u(x)\geqslant u(x^\circ)$ and from $ u(x^\circ)\geqslant u(x)$  we have $u(x) = u(x^\circ)$. It remains to see that $x^\circ$ is efficient. If 
$u(x^\prime)\geqslant u(x^\circ)$ then $u(x^\prime)\geqslant u(x)$ and, by the definition of $x^\circ$, 
$x^\prime\geqslant x^\circ$.\hfill$\Box$

\begin{defn}\label{defquasileon} The utility function $u: X\to\Lambda$ is quasi-Leontief if, for all $x\in X$,  there exists a point 
$x^\circ\in X$ such that 
\begin{equation}\label{equdefquasileon}
 \{x^\prime\in X : u(x^\prime)\geqslant u(x)\} =  \{x^\prime\in X : x^\prime\geqslant x^\circ\} 
 \end{equation}
 \end{defn}
 
 The defining relation (\ref{equdefquasileon}) in \ref{defquasileon} can also be written 
 \begin{equation}
u^{-1}\big(\uparrow\!\!u(x)\big) = \uparrow\!\!(x^\circ)
 \end{equation}\label{equdefquasileonbis}
 from which the uniqueness of $x^\circ$ can easily be inferred. We will denote by $u^0$ the map from 
 $X$ to $X$ defined by $u^0(x) = x^\circ$. Definition \ref{defquasileon} can be stated as follows 
 
 \medskip\noindent
 {\it $u: X\to\Lambda$ is quasi-Leontief if and only if there exists a function $u^\circ : X\to X$ such that, for all $x\in X$,}
 \begin{equation}\label{defquasileon}
 u^{-1}\big(\uparrow\!\!u(x)\big) = \uparrow\!\!\big(u^\circ(x)\big)
  \end{equation}

\begin{prop}\label{isoquasileon}
 Let $u: X\to\Lambda$ be a quasi-Leontief function then: \\
 $(1)$ $u$ is isotone and, for all $x\in X$, $x^\circ$ is the unique efficient point such that $u(x) = u(x^\circ)$ ; \\
 $(2a)$ $u^0$ is isotone and, for all $x\in X$, \\
 $(2b)$ $u^\circ(x)\leqslant x$; \\
 $(2c)$ $u^\circ\big(u^\circ(x)\big) = u^\circ(x)$; \\
 $(3)$ $u^0(X) = {\mathcal E}(u; X)$ and $x$ is efficient if and only if $u^0(x) = x$;  \footnote{$u^0$ a retraction from $X$ to ${\mathcal E}(u; X)$}\\  
 $(4)$ for all $x\in X$, $u^\circ(x) = \min\{x^\prime\in X: u(x^\prime)\geqslant u(x)\}$.
 
 \end{prop}
 
 {\it Proof:} Assume that $u: X\to\Lambda$ is quasi-Leontief and let, for all $x\in X$, $x^\circ$ be the      point such that (\ref{equdefquasileonbis}) holds.  From $u(x)\geqslant u(x)$ we have $x \geqslant x^\circ$. If $x_2\geqslant x_1$  then 
$x_2\geqslant x_1^\circ$  and therefore $u(x_2)\geqslant u(x_1)$; this proves the first part of $(1)$ and $(2b)$.   \\
From $x^\circ\geqslant x^\circ$ we have 
 $u(x^\circ)\geqslant u(x)$ and, since $u$ is isotone, we also have $u(x)\geqslant u(x^\circ)$ and therefore $u(x) = u(x^\circ)$. \\   If $u(x^\prime)\geqslant u(x^\circ)$ then $u(x^\prime)\geqslant u(x)$ and therefore $x^\prime\geqslant x^\circ$; we have shown that $x^\circ$ is efficient. Lemma \ref{efficient2} completes the proof of $(1)$. \\
 If  $u(x_2)\geqslant u(x_1)$ then $u\big(u^\circ(x_2)\big)\geqslant u(x_1)$ and therefore 
 $u^\circ(x_2)\geqslant u^\circ(x_1)$; this proves $(2a)$. To obtain $(2c)$ notice that  $u^\circ\big(u^\circ(x)\big)$ is efficient, that  $u\big(\, u^\circ\big(u^\circ(x)\big)\, \big) =  u\big(u^\circ(x)\big)$ and that 
 $u\big(u^\circ(x)\big) = u(x)$; therefore $u^\circ\big(u^\circ(x)\big)$ is efficient and 
 $u\big(\, u^\circ\big(u^\circ(x)\big)\, \big) =  u(x)$ which, from the second part of $(1)$, implies that 
 $u^\circ\big(u^\circ(x)\big) = u^\circ(x)$. $(3)$ very easily follows from $(1)$ and $(2)$ and $(4)$ follows  from 
 $(1)$ and Lemma \ref{efficient2}. \hfill$\Box$\footnote{A map from a partially ordered set to itself for which $(2a)$, $(2b)$ and $(2c)$ hold is an interior operator; the best known example of such an operator is the map which assigns to subsets $S$ of a given topological space $T$ their interior, often denoted by 
 $\stackrel{\circ}{S}$.   }

\begin{expl}[Examples of quasi-Leontief functions]\label{exquaileon1}\quad\hfill\\
$\boldsymbol{[1]}$({\bf Classical Leontief utility functions}) Take $X = {\mathbb R}^n$, $\Lambda = {\mathbb R}$ and let  $u(
x)=\min_{i\in[n]}a_ix_i$\, \footnote{$[n] = \{1, \cdots, n\}$} with $a = (a_1, \cdots, a_n)\in  {\mathbb R}^n_{++}$. Then 
$u(x^\prime)\geqslant u(x)$ if and only if, for all $j$, $a_jx^\prime_j\geqslant \min_{i\in[n]}a_ix_i$ that is 
$x_j\geqslant x^\circ_j$ with 
\begin{equation}
x^\circ_j = \frac{\min_{i\in[n]}a_ix_i}{a_j}
\end{equation}
$x\in\mathcal{E}(u, {\mathbb R}^n)$ if and only if 
\begin{equation}\label{careffex1}
\forall i, j\in [n]\quad a_ix_i = a_jx_j.
\end{equation}
$\boldsymbol{[2]}$ $X = {\mathbb R}^n_+$, $\Lambda$ and $a  = (a_1, \cdots, a_n)$ are as above and 
$u(
x)=\min_{i\in[n]}a_ix_i^{\alpha_i}$. Then 
\begin{equation}
x^\circ_j = \left[\frac{\min_{i\in[n]}a_ix_i^{\alpha_i}}{a_j}\right]^{1/\alpha_j}
\end{equation}
$x\in\mathcal{E}(u, {\mathbb R}^n_+)$ if and only if 
\begin{equation}\label{careffex2}
\forall i, j\in [n]\quad a_ix_i^{\alpha_i} = a_jx_j^{\alpha_j}.
\end{equation}

\noindent
 $\boldsymbol{[3]}$ $X$ is $\Real^n$ and 
 $p_1,p_2,...,p_n\in\mathbb{R}_+^n$ are $n$ linearly independent price vectors. Let
 $P=(p_{i,j})_{\substack{i=1,...,n \\j=1,...,n}}$ be the corresponding
 price matrix. Define the partial order $\geqslant_P$ on $\mathbb{R}^n$ by  $x\geqslant_Py$ if and only if $p_i\cdot x\geq p_i\cdot
 y$ for each $i=1,...,n$ and let  $u(x)=\min\{p_i\cdot
x: i\in [n]\}$. Let $\boldsymbol{1_n}$ be the element of $\mathbb{R}^n$ whose coordinates are all equal to $1$; $u(x^\prime)\geqslant u(x)$ can be written $Px^\prime\geqslant u(x)\boldsymbol{1_n}$. There exists a unique $x_P\in\mathbb{R}^n$ such that $Px_P = \boldsymbol{1_n}$; the inequality 
$u(x^\prime)\geqslant u(x)$ holds if and only if $x^\prime\geqslant_P u(x)x_P$ which shows that $u$ is quasi-leontief and that $u^\circ(x) = u(x)x_P$.

\noindent
 $\boldsymbol{[4]}$ If $u : X\to\mathbb{R}$ is a quasi-Leontief function then, for all $(a, b)\in \mathbb{R}_{++}\times\mathbb{R}$, $au + b$ is also a quasi-Leontief function. Indeed, $au(x^\prime) + b\geqslant au(x) + b$ is equivalent to $u(x^\prime)\geqslant u(x)$ and therefore $(au + b)^\circ = u^\circ$.
 
\end{expl}

\begin{prop}[Characterization of quasi-Leontief functions]\label{carregquasileon}\quad\\
Let $\Lambda_0$ be a subset of $\Lambda$ such that $u(X)\subset \Lambda_0\subset \downarrow\!\!\big(u(X)\big)$ and consider the following statements  : \\
$(1)$ there exists a function $u^\sharp : \Lambda_0 \to X$ such that 
\begin{equation}\label{equacarquasileon1}
\forall (x, \lambda)\, \in\,  X\, \times\Lambda_0\quad \left[x\geqslant u^\sharp(\lambda)\Leftrightarrow u(x)\geqslant \lambda\right]
\end{equation}
$(2)$ $u$ is isotone and there exists an isotone function $u^\sharp : \Lambda_0\to X$ such that
\begin{equation}\label{equacarquasileon2}
\forall (x, \lambda)\, \in\,  X\, \times\Lambda_0\quad \left[x\geqslant u^\sharp(u(x))\,\,\hbox{and}\,\, u(u^\sharp(\lambda))\geqslant \lambda\right]
\end{equation}
$(3)$ For all $\lambda\in\Lambda_0$ the set $\{x\in X: u(x)\geq\lambda\}$ has a smallest element,that is 
\begin{equation}\label{equacarquasileon3}
\forall \lambda\in\Lambda_0\quad \exists \lambda^\sharp\in X \hbox{ such that  } 
u^{-1}\{\uparrow\!\!\lambda\} = \uparrow\!\!(\lambda^\sharp)
\end{equation}

\medskip\noindent Then $(1)$, $(2)$ and $(3)$ are equivalent and they imply that  $u$ is quasi-Leontief; furthermore, if $\Lambda_0 = u(X)$ then they are all equivalent to $u$ being quasi-Leontief. \end{prop}

 {\it Proof:}
 Let us assume that $(1)$ holds. For all $x\in X$ let $u^\circ(x) = u^\sharp\big(u(x)\big)$. Then, for all 
 $(x^\prime, x)\in X\times X$, $u(x^\prime)\geqslant u(x)$ if and only if $x^\prime\geqslant u^\circ(x)$. That is, $u$ is quasi-Leontief and consequently isotone. Let us see that $u^\sharp$ is also isotone. \\ 
 Apply (\ref{equacarquasileon1}) with $x = u^\sharp(\lambda^\prime)$ to obtain, for all $\lambda\in  \Lambda_0$,  
 $u^\sharp(\lambda^\prime)\geqslant u^\sharp(\lambda)$ if and only if $u\big(u^\sharp(\lambda^\prime)\big)\geqslant \lambda$; taking $\lambda = \lambda^\prime$ yields, for all $\lambda^\prime$, $u\big(u^\sharp(\lambda^\prime)\big)\geqslant \lambda^\prime$. Consequently, if $\lambda^\prime\geq \lambda$ we have 
 $u\big(u^\sharp(\lambda^\prime)\big)\geqslant \lambda$ and therefore $u^\sharp(\lambda^\prime)\geqslant u^\sharp(\lambda)$.
 
 \medskip\noindent
 Let us see that (\ref{equacarquasileon2}) holds. We have just seen that, for all 
 $\lambda\in  \Lambda_0$,  $u\big(u^\sharp(\lambda)\big)\geqslant \lambda$. Taking $\lambda = u(x)$ in 
 (\ref{equacarquasileon1}) we obtain 
 $x\geqslant u^\sharp(u(x))$. We have shown that $(1)$ implies $(2)$.
 
 \medskip\noindent
 Let us assume that $(2)$ holds and that $u(x)\geqslant\lambda$. Since $u^\sharp$ is isotone we have $u^\sharp(u(x))\geqslant u^\sharp(\lambda)$ and therefore, by (\ref{equacarquasileon2}), 
 $x\geqslant u^\sharp(\lambda)$. One shows similarly that $x\geqslant u^\sharp(\lambda)$ implies $u(x)\geqslant\lambda$. We have shown that $(2)$ implies $(1)$.
 
 \medskip\noindent
 That $(1)$ and $(3)$ are equivalent is clear (which shows once more that $(1)$ implies that $u$ is quasi-Leontief, since $u(X)\subset\Lambda_0$). 
 
 \medskip\noindent
 Let us assume that $u:X\to\Lambda$ is quasi-Leontief and let $\Lambda_0 = u(X)$. For all 
 $\lambda\in u(X)$ we have $u\big(u^{-1}(\{\lambda\})\big) = \{\lambda\}$.  Since for all $x\in X$, 
 $x ^\circ = u^\circ(x)$ is the unique efficient point for which 
 $u(x^\circ) = u(x)$, $u^\circ$ is constant on $u^{-1}(\{\lambda\})$; let 
 $u^\sharp(\lambda)$ be the point of $X$ for which 
 $u^\circ\big(u^{-1}(\{\lambda\})\big) = \{u^\sharp(\lambda)\}$.  
  By construction we have $u\big(u^\sharp(\lambda)\big)= \lambda$ and $u^\sharp\big(u(x)\big) = u^\circ(x)$ and therefore  $x\geqslant u^\sharp\big(u(x)\big)$. We have shown that (\ref{equacarquasileon2}) holds. \hfill$\Box$
 
 \bigskip\noindent
 If $u: X\to\Lambda$ is a quasi-Leontief function then (\ref{equacarquasileon1}) and (\ref{equacarquasileon2}) hold with 
 $\Lambda_0 = u(X)$ and $\downarrow\!\!\big(u(X)\big)$ is the largest subset for which they could hold; the largest set for which they hold is 
 
 \begin{center}$\Lambda(u) = \{\lambda\in \downarrow\!\!\big(u(X)\big):  \exists \lambda^\sharp\in X \hbox{ such that  } 
u^{-1}\{\uparrow\!\!\lambda\} = \uparrow\!\!(\lambda^\sharp)\}$\end{center}.

 \begin{defn}\label{regularquasileon}
  A quasi-leontief function $u: X\to\Lambda$ is regular if, for all $\lambda\in\Lambda$, either 
 $u^{-1}(\uparrow\!\!\lambda) = \emptyset$ or $u^{-1}(\uparrow\!\!\lambda)$ has a smallest element.
 \end{defn}

  \noindent
 If $u: X\to\Lambda$ is a quasi-Leontief function then (\ref{equacarquasileon1}) and (\ref{equacarquasileon2}) with 
 $\Lambda_0 = u(X)$ and $\downarrow\!\!\big(u(X)\big)$ is the largest subset for which they could hold.
 
 \begin{lem}\label{trivregular } A quasi-Leontief function $u: X\to\Lambda$ is regular if and only if  (\ref{equacarquasileon1}), or  (\ref{equacarquasileon2}), holds with 
 $\Lambda_0 =  \downarrow\!\!\big(u(X)\big)$.
 \end{lem}
 
 {\it Proof:}  If (\ref{equacarquasileon1}) holds with $\Lambda_0 =  \downarrow\!\!\big(u(X)\big)$ then, for all $\lambda\in \downarrow\!\!\big(u(X)\big)$, $u^\sharp(\lambda)$ is the smallest element of 
 $u^{-1}(\uparrow\!\!\lambda)$; if $\lambda\not\in \downarrow\!\!\big(u(X)\big)$ then $u^{-1}(\uparrow\!\!\lambda) = \emptyset$. 
 
 \medskip\noindent Assume now that $u$ is a regular quasi-Leontief function. If $\lambda\in \downarrow\!\!\big(u(X)\big)$ then 
 $u^{-1}(\uparrow\!\!\lambda)$ is not empty; it has a smallest element, call it $\lambda^\sharp$. We trivially have $u^{-1}(\uparrow\!\!\lambda)\subset 
 \uparrow\!\!(\lambda^\sharp)$. If $x\geqslant\lambda^\sharp$ then 
 $u(x)\geqslant u(\lambda^\sharp)$, since $u$ is isotone, and 
 $u(\lambda^\sharp)\geqslant\lambda$, by definition of $\lambda^\sharp$; we have shown that $u(x)\geqslant\lambda$ and therefore that 
 $u^{-1}(\uparrow\!\!\lambda) = 
 \uparrow\!\!(\lambda^\sharp)$.\hfill$\Box$

 \begin{expl}[Non regular quasi-Leontief functions]\label{nonregularleon}\quad\hfill\\
$\boldsymbol{[1]}$ Let $X = [0, 1]\cup ]2, 3]$, $\Lambda = \mathbb{R}$ and 
$u(x) = x$. Obviously, $u$ is quasi-Leontief with $u^\circ = Id_X$ but 
$\{x\in X: u(x)\geqslant 2\} = ]2, 3]$ does not have a least element and 
$2\in\downarrow\!\!u\big(X\big)$. The trouble here is caused by the bounded decreasing sequence $2 + (1/n)$ which does not have an infimum in $X$. 

\medskip\noindent
$\boldsymbol{[2]}$ Let $X = [0, 1]\cup [2, 3]$, $\Lambda = \mathbb{R}$ and 
$u(x) = x$ for $x\in[0, 1]\cup ]2, 3]$ and $u(2) = 1.5$; $u$ is quasi-Leontief with $u^{\circ}(x) = x$ but  $\{x\in X: u(x)\geqslant 1.7\} = ]2, 3]$ does not have a least element and 
$1.7\in\downarrow\!\!u\big(X\big)$. The trouble here is due to the fact that $u$ is not upper semicontinuous. 
\end{expl}

 \medskip\noindent 
 The following proposition recapitulates and completes some of the properties of quasi-Leontief function. 
 
 \begin{prop}\label{recap}
 Let $u: X\to\Lambda$ be a quasi-Leontieff function and  let $\Lambda_0$ be a subset of $\Lambda$ such that $u(X)\subset \Lambda_0\subset \downarrow\!\!\big(u(X)\big)$.  
 
 \medskip\noindent
 $(1)$  $u^\circ = u^\sharp\circ u$\quad  and\quad    $\mathcal{E}(u; X) = \big\{ u^\sharp(\lambda) : \lambda\in\Lambda_0 \big\}$;  
 
  \medskip\noindent
 $(2)$  the function $\bar{u} = u\circ u^\sharp$ has the following properties:   it is   
   isotone and, for all 
  $\lambda\in \Lambda_0$, $\lambda\leqslant\bar{u}(\lambda)$ and 
 $\bar{u}\big(\bar{u}(\lambda)\big) = \bar{u}(\lambda)$; \footnote{$\bar{u}$ is a closure operator on $\downarrow\!\!\big(u(X)\big)$.}
 
  \medskip\noindent
 $(3)$ $u^\sharp$ is an order preserving bijection from $u(X) = \{\lambda : \bar{u}(\lambda) = \lambda\}$ to $\mathcal{E}(u; X)$ whose inverse is the restriction of $u$ to $\mathcal{E}(u; X)$.

 \medskip\noindent
 $(4)$ 
 $u(x) = \max\{\lambda : x\geqslant u^\sharp(\lambda)\}$ \quad and \quad  $u^\sharp(\lambda) = \min\{x : u(x)\geqslant\lambda\}$.

 \end{prop}
 
  {\it Proof:} The identity $u^\circ = u^\sharp\circ u$ is contained in the very first part of the proof of 
  Proposition \ref{carregquasileon}. If $x\in \mathcal{E}(u; X)$ if and only if $u^\circ(x) = x$, if and only if 
  $u^\sharp\big(u(x)\big) = x$; this shows that $\mathcal{E}(u; X) \subset \big\{ u^\sharp(\lambda) : \lambda\in \Lambda_0 \big\}$. The inclusion 
  $\mathcal{E}(u; X) \supset \big\{ u^\sharp(\lambda) : \lambda\in \Lambda_0 \big\}$ is a consequence of the identity $u^\sharp\circ u\circ u^\sharp= u^\sharp$ which can be proved as follows: 
  from $u^\circ(x)\leqslant x$ we have $u^\sharp\circ u\circ u^\sharp \leqslant u^\sharp$; from (\ref{equacarquasileon2})  we have $\lambda\leqslant u\big( u^\sharp(\lambda)\big)$ and since $u^\sharp$ is isotone we also have $u^\sharp(\lambda)\leqslant u^\sharp\big(u\big( u^\sharp(\lambda)\big) \big)$. This proves $(1)$.
  
  \medskip\noindent Both $u$ and $u^\sharp$ are isotone therefore $\bar{u}$ is isotone.  The relation 
  $\lambda\leqslant \bar{u}(\lambda)$ is part of (\ref{equacarquasileon2}) and, from   $u^\sharp\circ u\circ u^\sharp= u^\sharp$ we have  $u\circ u^\sharp\circ u\circ u^\sharp=u\circ u^\sharp$ that is 
  $\bar{u}\circ\bar{u} = \bar{u}$. This proves $(2)$. 
  
   \medskip\noindent We have seen that $\mathcal{E}(u; X) = \big\{ u^\sharp(\lambda) : \lambda\in \Lambda_0 \big\}$; the identity $\mathcal{E}(u; X) = \big\{ u^\sharp(\lambda) : \lambda\in u(X) \big\}$ follows from $u^\sharp(\lambda) = u^\sharp(\lambda^\prime)$ with 
   $\lambda^\prime  = u\big( u^\sharp(\lambda)\big)$. \\
   If $\bar{u}(\lambda) = \lambda$ then, from the definition of $\bar{u}$, we have $\lambda\in u(X)$.  If 
   $\lambda = u(x)$ then  $\bar{u}(\lambda) = u\big(u^\sharp\big(u(x)\big)\big) = \lambda$   since 
   $u\circ u^\sharp\circ u = u$;  notice that (\ref{equacarquasileon2}) and the fact that $u$ is isotone yield    $u\circ u^\sharp\circ u \leqslant u$ and that $u\leqslant u\circ u^\sharp\circ u$ follows from 
   $u\leqslant \bar{u}\circ u$.  We have shown that $u(X) = \{\lambda : \bar{u}(\lambda) = \lambda\}$.\\
   We have seen in the proof of Proposition \ref{carregquasileon} that   $u\big(u^\sharp(\lambda)\big)= \lambda$ and $u^\sharp\big(u(x)\big) = u^\circ(x)$; if $x\in \mathcal{E}(u; X)$ then $u^\circ(x) = x$. This proves $(3)$.
     
     \medskip\noindent $(4)$ follows directly from (\ref{equacarquasileon1}) and $(3)$. \hfill$\Box$
     
     \medskip\noindent In Proposition \ref{recap} one can always take 
     $\Lambda_0 = u\big(X\big)$ and $\downarrow\!\!u\big(X\big)$ if $u$ is regular. 
     
     The function $u^\sharp : u(X)\to X$ can be seen as a one to one isotone parametrization of $\mathcal{E}(u; X)$; for the standard Leontief functions,  this ``parametrized 
     path'' is the economist's {\bf expantion path} of $u$.
     
     \begin{expl}[More examples of quasi-Leontief functions]\label{exquaileon2}\quad\hfill\\
 $\boldsymbol{[1]}$ A classical Leontief utility function  $u(
x)=\min_{i\in[n]}a_ix_i$ defined on $\mathbb{R}^n$ with $a = (a_1, \cdots, a_n)\in  {\mathbb R}^n_{++}$ is regular with  $u^\sharp(\lambda) = (\lambda/a_1, \cdots, \lambda/a_n)$.

\medskip\noindent  Example $(2)$ from \ref{exquaileon1} is regular as well as $(3)$, with $u^\sharp(\lambda) = \lambda x_P$.

\medskip\noindent
$\boldsymbol{[2]}$ {\bf (Quasi-Leontief functions on a closed subset of $\boldsymbol{\mathbb{R}^n}$. )} \\If a function $u : X\to\mathbb{R}$ defined on a closed subset $X$ of  $\mathbb{R}^n$ is a regular quasi-Leontief function then \\\ $(1)$ it is increasing,\\
$(2)$ it has lower bounded upper level sets and\\
$(3)$ it is upper semicontinuous .

\medskip\noindent
For $n = 1$, conditions $(1)$, $(2)$ and $(3)$ imply that $u$ is a quasi-Leontief function.

\medskip\noindent
{\rm $(a)$ If $u$ is quasi-Leontief then it is increasing; if it is also regular then, for all 
$\lambda\in\mathbb{R}$, $\{x\in X: u(x)\geqslant\lambda\}$ is  a closed  subset of $X$ since it is either empty or it is $\{x\in X : x\geqslant 
u^\sharp(\lambda)\}$ and $u^\sharp(\lambda)$ is obviously a lower bound for  $\{x\in X: u(x)\geqslant\lambda\}$. \\
$(b)$ Assume that $n = 1$ and that $u$ is increasing and upper semicontinuous. \\The set 
$\{x\in X: u(x)\geqslant\lambda\}$ is closed in $\mathbb{R}$, since it is closed in $X$ and $X$ is closed in $\mathbb{R}$. If $\{x\in X: u(x)\geqslant\lambda\}$ has a lower bound it has a finite greatest lower bound   $\lambda^\sharp\in \{x\in X: u(x)\geqslant\lambda\}$ and since $u$ is increasing $\{x\in X: u(x)\geqslant\lambda\} = \{x\in X: x\geqslant\lambda^\sharp\}$.}

 \medskip\noindent{\rm  If $n > 1$ a characterization of regular quasi-Leontief functions is still possible. More generally, regular quasi-Leontief functions   on partially ordered sets are characterized by Theorem \ref{topleon1} below and by Proposition \ref{topleonsemi} when the partially ordered set is a  semilattice. }
 
 \medskip\noindent{\rm If $X$ is a convex subset of $\mathbb{R}^n$ then all quasi-Leontief functions $u : X\to\mathbb{R}$ are quasi-concave and upper-semicontinuous. But a partially ordered set, even if it a subset of a vector space, does not have to be convex, we will see that it does not even have to be a lattice.}

   \medskip\noindent
   $\boldsymbol{[3]}$ {\it  Let $u_i : X_i\to\Lambda$, $i\in [n]$, be a finite family of regular quasi-Leontief functions and let $\leqslant_i$ be the partial order on $X_i$ . Let $X = \prod_{i\in [n]}X_i$ be the product space endowed with the coordinatewise partial order and let 
   $u(x_1, \cdots, x_n) = \min_{i\in [n]}u_i(x_i)$}. Let us assume that $u^{-1}(\uparrow\!\!\lambda)$ is not empty;  
   then $(x_1, \cdots, x_n)\in u^{-1}(\uparrow\!\!\lambda)$ if and only if, for all 
   $i\in [n]$,  $u_i(x_i)\geqslant\lambda$ and therefore  
   $u^{-1}(\uparrow\!\!\lambda) = \prod_{i\in[n]}\!\uparrow\!\!u_i^\sharp(\lambda) = \uparrow\!\!(u_1^\sharp(\lambda), \cdots,  u_n^\sharp(\lambda))$. 
   In conclusion: \\ {\it $u$ is a regular quasi-Leontief function with  $u^\sharp(\lambda) = (u_1^\sharp(\lambda), \cdots, u_n^\sharp(\lambda))$\\ and \\$u^\circ(x_1, \cdots, x_n) = 
   (u_1^\sharp(\min_{i\in [n]}u_i(x_i)), \cdots, u_n^\sharp(\min_{i\in [n]}u_i(x_i)))$}.

    \medskip\noindent 
   $\boldsymbol{[4]}$ {\it Assume that arbitrary pairs of elements $(x_1, x_2)$ of $X$ always have  a least upper bound $x_1\vee x_2$}, for example, $X$ could be a lattice, and {\it let  $u_i  : X\to\Lambda$, $i\in [n]$, be regular quasi-Leontieff functions}. Let $u(x) = \min\{u_1(x), \cdots, u_n(x)\}$. Then $u(x)\geqslant \lambda$ if and only if, for all $i\in[n]$, 
   $u_i(x)\geqslant\lambda$ from which it is clear that either 
   $u^{-1}(\uparrow\!\!\lambda) = \emptyset$ or\\ 
   $u^{-1}(\uparrow\!\!\lambda) = \{x\in X: x\geqslant  u_1^\sharp(\lambda)\vee \cdots \vee u_n^\sharp(\lambda)\}$. In conclusion: \\   
  {\it  $u(x) = \min\{u_1(x), \cdots, u_n(x)\}$ is quasi-Leontief with\\   $u^\sharp(\lambda) = u_1^\sharp(\lambda)\vee \cdots \vee u_n^\sharp(\lambda)$ and\\ 
   $u^\circ(x) = u_1^\sharp(\min_{i\in[n]}u_i(x))\vee \cdots \vee u_n^\sharp(\min_{i\in[n]}u_i(x))$}.
   
    \medskip\noindent 
   $\boldsymbol{[5]}$ If $X$ has a smallest element $0_X$ then  constant maps   $u(x) = \lambda_0$ are  regular quasi-Leontief since $u^{-1}(\uparrow\!\!\lambda) \neq \emptyset$ if and only if $\lambda_0\geqslant \lambda$ in which case\\ 
   $u^{-1}(\uparrow\!\!\lambda) = X = \uparrow\!\!0_X$ with $u^\circ(x) = 0_X$ and, for $\lambda\geqslant\lambda_0$, $u^\sharp(\lambda) = 0_X$. \\
   If $X$ does not have a smallest element then $u(x) = \lambda_0$ can not be quasi-Leontieff since $u^{-1}(\uparrow\!\! u(x)) = X$.\\
  If $X$ is a lattice with smallest element $0_X$ then, from the conclusion of example $(3)$, $\min\{u(x), \lambda_0\}$ is a regular quasi-Leontieff function if $u : X\to\Lambda$ is regular quasi-Leontieff
 
  \medskip\noindent 
   $\boldsymbol{[6]}${\it  Let $u : \prod_{i\in [n]}X_i$ be a (regular) quasi-Leontief function then, for all $(\tilde{x}_2, \cdots \tilde{x}_n)\in \prod_{i > 1}X_i$, 
   $u_1(x_1) = u(x_1, \tilde{x}_2, \cdots \tilde{x}_n)$ is (regular) quasi-Leontief function on $X_1$}.
   
    \medskip\noindent
    {\rm     To see that $u_1$ is quasi-Leontief  let, for $x\in\prod_{i\in [n]}X_i$,  
    $u^\circ(x) = (u_1^\circ(x), \cdots, u_n(x))$ where $u_j^\circ(x)$ is the projection of $u^\circ(x)$ onto $X_j$ and notice that, from $x\geqslant u^\circ(x)$ we have, for all $x_1\in X_1$ and all $j > 1$, 
    $\tilde{x}_j\geqslant u_j^\circ(x_1, \tilde{x}_2, \cdots \tilde{x}_n)$ and therefore, $u_1(x_1^\prime)\geqslant u_1(x_1)$ if and only if 
    $x_1^\prime\geqslant u_j^\circ(x_1, \tilde{x}_2, \cdots \tilde{x}_n)$.
  
  \medskip\noindent  
    To see that $u_1$ is regular if $u$ is  notice that $\{x_1\in X_1: u_1(x_1)\geqslant\lambda\} = \emptyset$  if and only if either $\{x\in X: u(x)\geqslant\lambda\} = \emptyset$ or  
    $\{x\in X: u(x)\geqslant\lambda\} \neq \emptyset$ and $(x_1, \tilde{x}_2, \cdots \tilde{x}_n)\not\geqslant 
    u^\sharp(\lambda)$. \\
    Assume that $\{x_1\in X_1: u_1(x_1)\geqslant\lambda\} \neq \emptyset$ and let  $u^\sharp(\lambda) = (u_1^\sharp(\lambda), \cdots, u_n^\sharp(\lambda))$; taking an arbitrary point in $\{x_1\in X_1: u_1(x_1)\geqslant\lambda\}$ yields 
     $(\tilde{x}_2, \cdots \tilde{x}_n)\geqslant (u_2^\sharp(\lambda), \cdots, u_n^\sharp(\lambda))$ and therefore $u_1(x_1)\geqslant\lambda$ if and only if $x_1\geqslant u_1^\sharp(\lambda)$. }   
     
      \medskip\noindent 
   $\boldsymbol{[7]}${\it  If $u : X\to\Lambda$ is a (regular) quasi-Leontief function then, for all comprehensive subsets $S$ of $X$ the restriction 
   $u_{\vert S} : S\to\Lambda$ of $u$ to $S$ is a (regular) quasi-Leontief function}. 
   
   \medskip\noindent
   {\rm If $S$ is comprehensive then, for all $x\in S$, $u^\circ(x)\in S$; this shows that $u_{\vert S}$ is quasi-Leontief. Assume now that $u$ is regular and that $\{x\in S: u(x)\geqslant\lambda\}\neq\emptyset$; then, since 
   $\{x\in X: u(x)\geqslant\lambda\}\neq\emptyset$ and $S$ is comprehensive we have $u^\sharp(\lambda)\in S$ and consequently  
   $\{x\in S: u(x)\geqslant\lambda\} =\{x\in S: x\geqslant u^\sharp(\lambda)\}$. }
   \end{expl}

\subsection{A Characterization of Quasi-Leontief functions}\label{subsecarqf}
Assume that $u : X\to\Lambda$ is a quasi-Leontief function and let $x_1$ and $x_2$ be two arbitrary elements of $X$. Since 
$\mathcal{E}(u; X)$ is a totally ordered set, we can assume without loss of generality that  $u^\circ(x_1)\geqslant u^\circ(x_2)$ and, since $x_i\geqslant u^\circ(x_i)$ we also have $x_1\geqslant u^\circ(x_2)$. Since $u\big(u^\circ(x)\big) = u(x)$  we have $u\big(u^\circ(x_2)\big) = \min\{u(x_1), u(x_2)\}$. We have shown that an arbitrary quasi-Leontief function 
$u : X\to\Lambda$ has the following property: 

\medskip\noindent
{\bf (Property $\boldsymbol{\Phi}$)} {\it For all $x$ and $x^\prime$ in $X$ there exist  $x^{\prime\prime}\in X$ such that $x\geqslant x^{\prime\prime}$, 
$x^\prime\geqslant x^{\prime\prime}$ and $\min\{u(x), u(x^\prime)\} =  
u(x^{\prime\prime})$}.

\medskip\noindent Notice that the existence of a single quasi-Leontief function 
$u : X\to\Lambda$ on the partially ordered set $X$ implies that $X$ is a {\bf filtered partially ordered set}. \, \footnote{A  partially ordered set $L$ is filtered if, $\forall l_1, l_1\in L $ $\exists l_3\in L$ such that $l_1\geqslant l_3$ and $l_2\geqslant l_3$; property $\Phi$ implies that both $X$ and the graph of $u$  are filtered. }

\medskip\noindent
Assume that $\mathcal{C}$ is a chain in $X$ for which $\inf\mathcal{C}$ exists, call it $x_{\mathcal{C}}$. Since $u$ is monotone we have, for all 
$x\in \mathcal{C}$, $u(x)\geqslant u(x_{\mathcal{C}})$. Let $\lambda\in\Lambda$ such that, for all $x\in\mathcal{C}$, $u(x)\geqslant \lambda$. Assuming that $u$ is quasi-Leontief and regular we must have, for all $x\in\mathcal{C}$, 
$x\geqslant u^\sharp(\lambda)$ and consequently 
$x_{\mathcal{C}}\geqslant u^\sharp(\lambda)$. Taking the value on both sides yields 
$u(x_{\mathcal{C}})\geqslant u\big(u^\sharp(\lambda) \big)$ and therefore $u(x_{\mathcal{C}})\geqslant \lambda$, from $(2)$ of Proposition \ref{recap}.

\medskip\noindent
We have shown that an arbitrary regular quasi-Leontief function $u : X\to\Lambda$ has the following property:

\medskip\noindent
{\bf Property (CIP)} {\it If  $\mathcal{C}$ is a chain in $X$ for which $\inf\mathcal{C}$ exists then $u\big(\mathcal C\big)$ has a greatest lower bound in $\Lambda$ and 
$u(\inf\mathcal{C}) = \inf u\big(\mathcal{C}\big)$.} \footnote{CIP stands for ``chain inf preserving''.}

\medskip\noindent It is clear that a function for which (CIP) holds is isotone.

\medskip\noindent If $u$ is a regular quasi-Leontieff function then an arbitrary non empty upper level set 
$u^{-1}\big(\uparrow\!\!(\lambda)\big)$ is bounded below by $u^\sharp(\lambda)$. The three properties that we have listed above, $(\Phi)$, (CIP) and having lower bounded upper level sets, essentially characterize quasi-Leontief function. More precisely: 

\begin{prop}\label{algcarquasileon} Let $X$ be a filtered partially ordered set in which arbitrary chains that are bounded below have a greatest lower bound. Let $\Lambda$ be a totally ordered set. Then $u : X\to\Lambda$  is a regular quasi-Leontief function function if and only if it has properties $(\Phi)$, (CIP) and has, possibly empty, lower bounded upper level sets.
\end{prop}

 {\it Proof:}  We have already seen that a regular quasi-Leontief function has the three properties in question. Let $u : X\to\Lambda$ be a function  for which these properties hold and fix an arbitrary $\lambda\in\Lambda$. By hypothesis the set $u^{-1}\big(\uparrow\!\!(\lambda)\big)$  is bouded below consequently, if 
 ${\mathcal C}$ is an arbitrary chain in  $u^{-1}\big(\uparrow\!\!(\lambda)\big)$ then  
 $\inf{\mathcal C}$ exists in $X$; by (CIP) $\inf{\mathcal C}\in u^{-1}\big(\uparrow\!\!(\lambda)\big)$. By Zorn's Lemma the set of minimal elements of  $u^{-1}\big(\uparrow\!\!(\lambda)\big)$ is not empty. This set is of cardinality one; indeed, if $x$ and $x^\prime$ are minimal elements of $u^{-1}\big(\uparrow\!\!(\lambda)\big)$ we can, by Property $(\Phi)$ find 
 $x^{\prime\prime}\in\downarrow\!\!\{x, x^\prime\}$ such that 
 $u(x^{\prime\prime}) = \min\{u(x), u(x^\prime)\}$; since both $x$ and $x^\prime$ are in $u^{-1}\big(\uparrow\!\!(\lambda)\big)$ we also have 
 $u(x^{\prime\prime})\geqslant\lambda$. By minimality we must have 
 $x = x^{\prime\prime}$ and $x^\prime = x^{\prime\prime}$. Let $\lambda^\sharp$ be the unique minimal element of $u^{-1}\big(\uparrow\!\!(\lambda)\big)$. To complete the proof we have to see that  $u^{-1}\big(\uparrow\!\!(\lambda)\big) = \uparrow\!\!(\lambda^\sharp)$. By definition of $\lambda^\sharp$ we have 
 $u^{-1}\big(\uparrow\!\!(\lambda)\big) \subset \uparrow\!\!(\lambda^\sharp)$.
 
 Let $x$ be an arbitrary element of $\uparrow\!\!(\lambda^\sharp)$; by 
 (CIP), which implies that $u$ is isotone, we have $u(x)\geqslant u(\lambda^\sharp)$. From $\lambda^\sharp\in u^{-1}\big(\uparrow\!\!(\lambda)\big) $ we have  $u(\lambda^\sharp)\geq\lambda$ and consequently $u(x)\geq\lambda$.\hfill$\Box$
 
 \begin{lem}\label{minchaine} Let $X$ be a filtered partially ordered set endowed with a topology for which intervals are compact. If for all $x\in X$ the set  $\uparrow\!\!(x)$ is closed then all   chains in $X$ which are bounded below  have a greatest lower bound. Furthermore, the greatest lower bound will belong to all closed subsets of $X$ in which the chain is contained.
  \end{lem}
 
  {\it Proof:} Let $\mathcal{C}$ be a chain in $X$ and let $\underline{x}$ be a lower bound of $\mathcal{C}$. For all $x\in \mathcal{C}$ the  set  $\cap_{x\in\mathcal{C}}[x, \underline{x}]$ is  compact and not empty, since it contains $\underline{x}$; call it $C(\underline{x})$. For all $y\in C(\underline{x})$ the set 
  $\uparrow\!\!(y)\cap C(\underline{x})$ is compact and not empty; if  $\{y_0, \cdots, y_m\}$ are arbitrary elements of $C(\underline{x})$ then, for all $x\in\mathcal{C}$ the set $\cap_{i=0}^m[x, y_i]$ is compact and not empty, since it contains $x$. If $\{x_0, \cdots, x_n\}$ are elements of 
  $\mathcal{C}$ then  one of them is the smallest, since $\mathcal{C}$ is a chain; let $x_0 = \min\{x_0, \cdots, x_n\}$. From 
  $[x_i, y_j]\supset [x_0, y_j]$ we have 
  $\cap_{j=0}^n\cap_{i=0}^m[x_j, y_i] = \cap_{i=0}^m[x_0, y_i]$. This shows that the family of compact sets 
  $\{\cap_{i=0}^m[x, y_i]: x\in\mathcal{C}\}$ has the finite intersection property and consequently that 
  $\cap_{x\in\mathcal{C}}\cap_{i=0}^m[x, y_i]\neq\emptyset$. If 
  $x^\star$ is an arbitrary element of that intersection then 
  $x^\star\in \cap_{i=0}^m\big(\uparrow\!\!(y_i)\cap C(\underline{x})\big)$. We have shown that the family of compact sets $\{\uparrow\!\!(y)\cap C(\underline{x}): y\in C(\underline{x})\}$ has the finite intersection property; the set $\cap_{y\in C(\underline{x})}\big[\uparrow\!\!(y)\cap C(\underline{x})\big]$ is therefore not empty. Let $y^\star$ be an element of that intersection. Let us see that 
  $\cap_{y\in C(\underline{x})}\big[\uparrow\!\!(y)\cap C(\underline{x})\big] = \{y^\star\}$. If $\hat{y}\in\cap_{y\in C(\underline{x})}\big[\uparrow\!\!(y)\cap C(\underline{x})\big]$ then $\hat{y}\in C(\underline{x})$ and therefore $y^\star\geqslant \hat{y}$; interchanging $y^\star$ and $\hat{y}$ yields $\hat{y}\geqslant y^\star$. 
  
  \medskip\noindent To complete the first part of the proof let us see that  $y^\star$ is the greatest lower bound of $\mathcal{C}$. If $x^\prime$ is an arbitrary lower bound of $\mathcal{C}$ choose a point $x^{\prime\prime}\in\big[\downarrow\!\!(y^\star)\cap\downarrow\!\!(x^\prime)\big]$ and let $y^{\star\star}$ be an element of $\cap_{y\in C(x^{\prime\prime})}\big[\uparrow\!\!(y)\cap C(x^{\prime\prime})\big]$. Since  $y^\star$ is a lower bound of $\mathcal{C}$ and $y^\star\geqslant x^{\prime\prime} $ we have $y^\star\in C(x^{\prime\prime})$ and therefore $y^{\star\star}\geqslant y^\star$; this imlies that $y^{\star\star}$ is a lower bound of 
  $\mathcal{C}$ such that $y^{\star\star}\geqslant \underline{x}$, in other words, $y^{\star\star}\in \cap_{y\in C(\underline{x})}\big[\uparrow\!\!(y)\cap C(\underline{x})\big]$ and therefore $y^{\star\star} = y^\star$.\\
  By construction we have $x^\prime\geqslant x^{\prime\prime}$ and by hypothesis $x^\prime$ is a lower bound of $\mathcal{C}$, that is 
  $x^\prime\in C(x^{\prime\prime})$ and therefore 
  $y^{\star\star}\geqslant x^\prime$. We have shown that 
  $y^{\star}\geqslant x^\prime$. 
  
  \medskip\noindent Let now $F$ be a closed subset of $X$ such that 
  $\mathcal{C}\subset F$. For all $x\in \mathcal{C}$ the set 
  $[x, y^\star]\cap F$ is compact and not empty. As above one shows that the family $\{[x, y^\star]\cap F: x\in  \mathcal{C}\}$ has the finite intersection property; by compactness we can find an element $y^{\star\star}$ in $\cap_{x\in  \mathcal{C}}[x, y^\star]\cap F$ and this implies that $y^{\star\star}$ is a lower bound of $\mathcal{C}$ such that $y^{\star\star}\geqslant y^\star$ and therefore that  $y^{\star\star} = y^\star$ and finally that $y^\star\in F$.  \hfill$\Box$
  
  \bigskip
  Let us say that a partial order on a topological space  is an {\bf upper semi continuous partial order} if all set of the form $\uparrow\!\!(x)$ are closed.
  
  \begin{thm}\label{topleon1}  Let $X$ be a filtered partially ordered set.   Assume that both $X$ and  $\Lambda$ are endowed with a topology for which intervals are compact and the partial orders are upper semicontinuous. Let $u : X\to\Lambda$ be a function with closed and lower bounded upper level sets. Then $u$ is a regular quasi-Leontief function if and only is it is isotone and  it has property $(\Phi)$.
  \end{thm}
  
  {\it Proof:} We establish the non obvious part of the theorem. We have to show that (CIP) holds. By Lemma \ref{minchaine} chains in $X$ which are bounded below have a greatest lower bound; let $\mathcal{C}$ be such a chain.  For all 
  $x\in \mathcal{C}$ we have $u(x)\geqslant u(\inf\mathcal{C})$ and therefore $\inf\!u\big(\mathcal{C}\big)$ exists in $\Lambda$ and 
  $\inf\!u\big(\mathcal{C}\big)\geqslant u(\inf\mathcal{C})$. The chain 
  $\mathcal{C}$ is contained in $\{x\in X: u(x) \geqslant \inf\!u\big(\mathcal{C}\big)\}$ and this set is closed, which implies that 
  $\inf\mathcal{C}$ belongs to $\{x\in X: u(x) \geqslant \inf\!u\big(\mathcal{C}\big)\}$; we have shown that   
  $u(\inf\mathcal{C})\geqslant \inf\!u\big(\mathcal{C}\big)$. \hfill$\Box$

\subsection{Quasi-Leontief functions on semilattices}\label{subsecsemilat}
An {\bf inf-semilattice} is a partially ordered set $(L, \geqslant)$ for which sets of cardinality two $\{x_1, x_2\}$ allways have a greatest lower bound, written   $x_1\wedge x_2$. A totally ordered set is an inf-semilattice; $\mathbb{R}^n$ and    $\mathbb{R}^n_+$ are inf-semilattices;   $\{(x_1, x_2)\in \mathbb{R}^2: x_1 + x_2\leqslant 1\}$ is an inf-semilattice which is not a sublattice of $\mathbb{R}^n$. A {\bf topological inf-semilattice} is a semilattice $L$ equipped with a topology for  for which the map $(x_1, x_2)\mapsto x_1\wedge x_2$ is continuous. The examples given above are all instances of topological inf-semilattices. Let $X$ be an inf-semilatice; for all 
$x\in X$ we have $\uparrow\!\!(x) = \{y\in X: x\wedge y = x\}$ and  
$\downarrow\!\!(x) = \{y\in X: x\wedge y = y\}$; if $X$ is a topological inf-semilattice then, from the continuity of the map $y\mapsto x\wedge y$, we have that both $\uparrow\!\!(x)$ and $\downarrow\!\!(x)$ are closed. Since  
$\Lambda$ is totally ordered it is also an inf-semilattice; notice that a function $u : X\to\Lambda$ such that, for all $x_1, x_2\in X$, 
 $u(x_1\wedge x_2) = \min\{u(x_1), u(x_2)\}$ is isotone. Such a function is a {\bf inf-semilattice homomorphism}.
 
 \begin{lem}\label{inftreilmor} A quasi-Leontief function $u : X\to\Lambda$ defined on an inf-semilattice $X$ is always an inf-semilattice homomorphism.
 \end{lem}
 
 {\it Proof:} Since a quasi-Leontief function has Property $(\Phi)$ let $x_1$ and $x_2$ be two arbitrary points of $X$ and let $x^{\prime}\in X$ such that 
 $x_1\geqslant x^{\prime}$, $x_2\geqslant x^{\prime}$ and\\ $u(x^{\prime}) = \min\{u(x_1), u(x_2)\}$. \\
 Since $u$ is isotone and $x_i\geqslant x_1\wedge x_2$ we have 
 $ \min\{u(x_1), u(x_2)\}\geqslant u(x_1\wedge x_2)$ and therefore 
 $u(x^\prime)\geqslant u(x_1\wedge x_2)$. \\
 From $x_i\geqslant x^\prime$ we have $x_1\wedge x_2\geqslant x^{\prime}$ and consequently $u(x_1\wedge x_2)\geqslant u(x^{\prime})$. We have shown that $u(x_1\wedge x_2) = u(x^{\prime})$, that is 
$u(x_1\wedge x_2) = \min\{u(x_1), u(x_2)\}$.

  \hfill$\Box$
  
  \begin{lem}\label{effsublat} If   $u : X\to\Lambda$ is a quasi-Leontief function defined on an inf-semilattice $X$ then $\mathcal{E}(u, X)$ is a sub-semilattice of $X$.
 \end{lem}
 {\it Proof:} By Lemma \ref{efficientchain}.
 \hfill$\Box$

 \bigskip\noindent From Theorem 
 \ref{topleon1} and Lemma \ref{inftreilmor}  we have the following characterization of regular quasi-Leontief functions on topological semilattices.
 
 \begin{prop}\label{topleonsemi}  Assume  that both $X$ and $\Lambda$ are  topological inf-semilattices for which intervals are compact.  Let $u : X\to\Lambda$ be a function with closed and lower bounded upper level sets. Then $u$ is a regular  quasi-Leontief function if and only if, for all $x_1, x_2\in X$, 
 $u(x_1\wedge x_2) = \min\{u(x_1), u(x_2)\}$.  \end{prop}
 
  \begin{cor}\label{carleonsubRn}
Let $X$ be a closed inf-semilattice of $\mathbb{R}^n$ and let $\Lambda$ be a closed subset of $\mathbb{R}$. Let $u : X\to\Lambda$ be an upper semicontinuous function with lower bounded upper level sets. Then $u$ is a regular quasi-Leontief if and only if, for all $x_1, x_2\in X$, $u(x_1\wedge x_2) = \min\{u(x_1), u(x_2)\}$. 
\end{cor}

\bigskip\noindent 
We have seen in $(3)$ of Examples \ref{exquaileon2}  that for all finite family $u_i : X_i\to\Lambda$, $i\in [n]$, of regular quasi-Leontief functions  the function $u(x_1, \cdots, x_n) = \min_{i\in [n]}u_i(x_i)$ is a regular quasi-Leontief function  on the product space $X = \prod_{i\in [n]}X_i$  endowed with the coordinatewise partial order. If all the $X_i$ are inf-semilattices then 
all the $u_i$ are semilattices homomorphisms. The next result shows that, under suitable but mild assumptions this is the generic case. 

\begin{prop}\label{minprodquasileon} Let $u :  \prod_{i\in [n]}X_i\to\Lambda$ be a (regular) quasi-Leontief function defined on a finite product of semilattices. Let $S\subset \prod_{i\in [n]}X_i$ be a   subset which has an upper bound $\bar{x}$ in $\prod_{i\in [n]}X_i$. Then, there exist a family $u_i : X_i\to\Lambda$, $i\in [n]$ of (regular) quasi-Leontief functions such that, for all $(x_1, \cdots, x_n)\in S$, 
$u(x_1, \cdots, x_n) = \min_{i\in [n]}u_i(x_i)$.
\end{prop}

{\it Proof:} Let $\bar{x}\in X$ be an upper bound of $S$. For all 
$x = (x_1, \cdots, x_n)\in \prod_{i\in [n]}X_i$ let  
$x_{[1]} = (x_1, \bar{x}_2 , \cdots, \bar{x}_n)$, 
  $x_{[2]} = (\bar{x}_1, x_2, \bar{x}_3,   \cdots, \bar{x}_n)$ and so on and notice that, for all $x\in S$,  $x = x_{[1]}\wedge x_{[2]}\wedge\cdots x_{[n]}$. Since $u$ is an inf-semilattice homomorphism we have, for all $x\in S$, 
  $u(x) =  \min\{u(x_{[1]}), \cdots, u(x_{[n]})\}$. For all $x\in  \prod_{i\in [n]}X_i$ let $u_1(x_1) = u(x_1, \bar{x}_2 , \cdots, \bar{x}_n)$, 
  $u_2(x_2) = u(\bar{x}_1, x_2,  \cdots, \bar{x}_n)$ and so on. We know from 
  $(6)$ of Examples \ref{exquaileon2} that $u_i$ is a (regular) quasi-Leontief 
  function on $X_i$.  \hfill$\Box$
  
  \begin{cor}\label{corminprodquasileon} Let $X_i\subset\mathbb{R}$, $i\in [n]$, be a finite family of closed intervals of $\mathbb{R}$ and let 
  $u : \prod_{i\in [n]}X_i\to \mathbb{R}$ be a regular quasi-Leontief function. 
  If $S\subset  \prod_{i\in [n]}X_i$ is a  subset of the product  with an upper bound in $\prod_{i\in [n]}X_i$ then there exists increasing upper semicontinuous functions 
  $u_i : X_i\to\mathbb{R}$ with lower bounded upper level sets such that, for 
  all $(x_1, \cdots, x_n)\in S$, 
  $u(x_1, \cdots, x_n) =  \min_{i\in [n]}u_i(x_i) $. 
  \end{cor}
  
  {\it Proof:} From Proposition \ref{minprodquasileon} and $(2)$ of Examples 
  \ref{exquaileon2}.\hfill$\Box$ 
  
  \bigskip
  We close this section with a characterization of Leontief functions. 
  
  \begin{prop}\label{homoquasileon} A regular  
  quasi-Leontief function $u : \mathbb{R}^n_{++}\to \mathbb{R}_{+}$ is a Leontief function if and only if it is homogeneous. 
  \end{prop}
   {\it Proof:} Let $\tilde{u}$ be the restriction of $u$ to the cube 
   $S = \{x\in \mathbb{R}^n_{++}: \max_{i\in[n]}x_i\leq 1\}$. By Corollary \ref{corminprodquasileon} we have, for all $(x_1, \cdots, x_n)\in S$, 
  $\tilde{u}(x_1, \cdots, x_n) =  \min_{i\in [n]}u_i(x_i) $ where each  function 
  $u_i :\,  ]0, 1]\to\mathbb{R}$ is increasing and upper semicontinous. Let us see that, for all $i\in [n]$, $u_i$ is homogeneous, that is $u_i(\lambda x_i) = \lambda u_i(x_i)$ if $0 < x_i \leqslant 1$ and $0 < \lambda \leqslant 1$. If $n=1$ there is nothing to prove since $u$ is homogeneous. We complete the proof by induction. 
  
  \medskip\noindent 
  Let $n = m+1$  and assume that for all homogeneous functions 
  $v : ]0, 1]^{m}\to\mathbb{R}_+$ of the form  $v(x_1, \cdots, x_m) = \min_{i\in[m]}v_i(x_i)$ the functions $v_i$ are homogeneous.

\medskip\noindent  
Define $v$ on $]0, 1]^{m}$ by 
   \begin{equation*}v(x_1, \cdots, x_{m}) = u(x_1, \cdots, x_{m}, \max_{i\in [m]}x_i).
  \end{equation*}
  Let us write $x_{[m]}$ for 
  $(x_1, \cdots, x_{m})$.   For $0 < \lambda \leqslant 1$ we have 
  \begin{eqnarray*}
v(\lambda x_{[m]})& = &u(\lambda x_{[m]}, \lambda\max_{i\in [m]} x_i)\\ 
&=& \lambda u(x_{[m]}, \max_{i\in [m]}x_i)\\
& = &
  \lambda v(x_{[m]}).
  \end{eqnarray*}
  
  From the  proof of Proposition \ref{minprodquasileon} we have 
  $v(x_1, \cdots, x_{n-1}) = \min_{i\in [n-1]}v_i(x_i)$ with 
  $v_i(x_i) = v({\bar x_{[n-1], i}})$ where ${\bar x_{[n-1], i}}$ has all its coordinates equal to $1$ with the exception of coordinate $i$ which is $x_i$. We have 
  $v({\bar x_{[n-1],i}}) = u({\bar x_{[n-1]}, i}, 1)$ and therefore, still from the 
  proof of Proposition \ref{minprodquasileon}, $v_i(x_i) = u_i(x_i)$. \\
  We have shown that, for $i\in [m]$, $u_i$ is homogeneous.    A permutation of the indices shows that $u_{m+1}$ is also homogeneous.  
  
  \medskip\noindent
  Let $a_i = u_i(1)$;  for $x\in ]0, 1]^n$ we have $u_i(x_i) = x_iu_i(1)$ and therefore\\  
  $u(x_1, \cdots, x_n) = \min_{i\in [n]}a_ix_i$. For an arbitrary $x\in\mathbb{R}^n_{++}$ choose $\lambda \leqslant 1$ such that   $\lambda x\in ]0, 1]^n$; we then have $\lambda u(x) = u(\lambda x) =   \min_{i\in [n]}a_i(\lambda x_i) = \lambda  \min_{i\in [n]}a_ix_i$ and finally\\
  $u(x) = \min_{i\in [n]}a_ix_i$. \hfill$\Box$

\subsection{Maximization of quasi-Leontief functions}\label{subsecmaxql}
If $u : X\to\Lambda$  is a quasi-Leontief function and if $S$ is an arbitrary subset of $X$ with a largest element $\bar{x}$ then 
$u(\bar{x}) = \max_{x\in S}u(x)$; also, $u^\circ(\bar{x})$ is the largest element of   $S\cap\mathcal{E}(u; X)$ and $u(\bar{x}) = u\big(u^\circ(\bar{x})\big)$. If $S$ does not have a largest element but if  $S\cap\mathcal{E}(u; X)$ has a largest element then $u$ still achieves its maximum value on $S$ as the following proposition shows. 

\begin{prop}\label{argmax1} Given a quasi-Leontief function $u : X\to\Lambda$ and a comprehensive subset $S\subset X$ the set 
$\arg\!\max(u; S)$ is not empty if and only if $S\cap\mathcal{E}(u; X)$ contains a largest element. Furthermore, if $\bar{x}$ is the largest element of $S\cap\mathcal{E}(u; X)$ then $\bar{x}\in\arg\!\max(u; S)$ and, for all 
$x\in\arg\!\max(u; S)$, $x\geqslant \bar{x}$ .
\end{prop}

 {\it Proof:} Assume that the set $\arg\!\max(u; S)$ is not empty and let $x^\star$ be an arbitrary element of $\arg\!\max(u; S)$. Since $S$ is comprehensive we have $u^\circ(x^\star)\in S$ and since 
 $u(x^\star)=  u\big(u^\circ(x^\star)\big) $ we also have 
 $u^\circ(x^\star)\in \mathcal{E}(u; X)$. If $x^\circ$ is an arbitrary element of $S$ we have $u\big(u^\circ(x^\star)\big)\geqslant u(x)$ and therefore $u^\circ(x^\star)\geqslant u^\sharp\big(u(x)\big)$. If 
 $x\in \mathcal{E}(u; X)$ then $u^\sharp\big(u(x)\big) = x$  which proves that $u^\circ(x^\star)$ is the largest element of 
 $S\cap\mathcal{E}(u; X)$.
 
 \medskip\noindent
 Assume now that $S\cap\mathcal{E}(u; X)$ is not emptyset and that it has a largest element $\bar{x}$. For all $x\in S$ we have 
 $u^\circ(x)\in S$ and therefore $\bar{x}\geqslant u^\circ(x)$; since $u$ is isotone and $u\big( u^\circ(x)\big) = u(x)$ we have 
 $u(\bar{x})\geqslant u(x)$. We have shown that $\bar{x}\in\arg\!\max(u; S)$. If $x$ is another element of $\arg\!\max(u; S)$ then $u(x)\geqslant u(\bar{x})$ and therefore $u \geqslant u^\sharp\big(u(\bar{x})\big) = u^\circ(\bar{x}) = \bar{x}$.  \hfill$\Box$ 
 
 \bigskip\noindent{\bf Remark 2.} If $S$ is a non empty comprehensive set then $u^\circ$ maps $S$ to itself. The set of fixed points of $u^0$, as a map from $S$ to $S$, is exactly $S\cap\mathcal{E}(u; X)$; Proposition \ref{argmax1} can be parphrased as follows: {\it $\arg\!\max(u; S)\neq\emptyset$ if and only if $u^\circ_{\vert S}$ has a largest fixed point.}
 
  \bigskip\noindent The next proposition completes Proposition \ref{argmax1}. 
  
  \begin{prop}\label{argmax2} Let $u : X\to\Lambda$ be a quasi-Leontief function and $S$ an arbitrary but non empty subset of $X$. Then, 
  
  \medskip\noindent
  $(1)$ $\arg\!\max(u; S)\neq\emptyset$ if and only if  
  $\arg\!\max(u; \downarrow\!\!(S))\neq\emptyset$.
  
   \medskip\noindent 
   $(2)$ If $S$ is comprehensive then the following assertions are equivalent:\\
   $(a)$ $\arg\!\max(u; S)\neq\emptyset$;\\
   $(b)$ there exists $x_0\in S$ such that 
   $\arg\!\max(u; \uparrow\!\!(x_0)\cap S)\neq\emptyset$;\\
   $(c)$ for all  $x\in S$  
   $\arg\!\max\big(u; \uparrow\!\!\big(u^\circ(x)\big)\cap S\big)\neq\emptyset$.
 \end{prop} 
 
 {\it Proof:} $(1)$ If $\bar{u}\in \arg\!\max(u; S)$ then 
 $\bar{u}\in\arg\!\max(u; \downarrow\!\!(S))$ since $u$ is isotone and for all $x\in \downarrow\!\!(S)$ there exists $\bar{x}\in S$ such that $\bar{x}\geqslant x$. Reciprocally, if $\bar{u}\in\arg\!\max(u; \downarrow\!\!(S))$ then choose $\bar{x}\in S$ such that $\bar{x}\geq \bar{u}$; since $u$ is isotone we have $\bar{x}\in\arg\!\max(u; \downarrow\!\!(S))\cap S$ and {\it a fortiori}  $\bar{x}\in \arg\!\max(u; S)$.  
 
 \medskip\noindent
 $(2)$ Assume that $S =\downarrow\!\!(S)$. \\ To see that $(a)$ implies $(b)$ notice that if $x$ is the largest  element of  $S\cap\mathcal{E}(u; X)$ then   
 $x\in\arg\!\max\big(u; \uparrow\!\!\big(u^\circ(x)\big)\cap S\big)$. \\
 We show that $(b)$ implies $(a)$. Let $\bar{x}$ be the largest element of $\big[\uparrow\!\!(x_0)\cap S\big]\cap\mathcal{E}(u; X)$. If $x^\star\in S\cap\mathcal{E}(u; X)$ then either $x^\star\geqslant \bar{x}$ or 
 $\bar{x}\geqslant x^\star$. If $x^\star\geqslant \bar{x}$ then 
 $x^\star\in \big[\uparrow\!\!(x_0)\cap S\big]\cap\mathcal{E}(u; X)$, since 
 $\bar{x}\geqslant x_0$, and consequently $x^\star = \bar{x}$. We have shown that $\bar{x}$ is the largest element of $S\cap\mathcal{E}(u; X)$; $(a)$ follows from  Proposition \ref{argmax1}.\\
 Let us see that $(a)$ implies $(c)$. If $S$ has a largest efficient $\bar{x}$ then, for all $x\in S$,  $\bar{x}\uparrow\!\!\big(u^\circ(x)\big)\cap S$ and $\bar{x}$  is the largest efficient point of $\uparrow\!\!\big(u^\circ(x)\big)\cap S$ . \\
 To conclude, $(c)$ implies $(b)$ trivially. \hfill$\Box$

 \bigskip\noindent
 {\bf Property (CUC)} {\it We will say that a subset $S$ of $X$ is  chain upper closed if all chains $\cal{C}$ in $S$ that have an upper bound in $X$ have a least upper bound $\sup\cal{C}$ and that least upper bound belongs to S.}
 
 \bigskip\noindent If $S$ is an arbitrary set  with a largest element $\bar{x}$ and if $u : S\to\Lambda$ is an arbitrary isotone function  then, trivially, $\bar{x}\in  \arg\!\max(u; S)$; if $S$ does not have a largest element but has a maximal element the same conclusion holds if $u$ is a quasi-Leontief function. Recall that $\bar{x}\in S$ is a maximal element if 
 $S\cap\uparrow\!\!\bar{x} = \{\bar{x}\}$. Let $\boldsymbol{Max}(S)$ be the, possibly empty, set of maximal elements of $S$
 
 \medskip\noindent The subset $\{(x_1, x_2)\in\mathbb{R}^2: 0 < x_1 
 \hbox{ and } x_1 + x_2\leqslant 1\}$ is CUC, it is bounded above and it does not have a largest element but is has plenty of maximal elements.
 
 \begin{thm}\label{maxcuc} Let $u : X\to\Lambda$ be a quasi-Leontief function. Then, for all comprehensive non empty (CUC) and bounded above subset $S$ of 
 $X$ 
  \begin{equation*}
 \boldsymbol{Max}(S)\, \cap\, \arg\!\max(u; S) \, \neq \, \emptyset
 \end{equation*}
 and consequently $ \max_{x\in S}u(x) =  \max_{x\in \boldsymbol{Max}(S)}u(x)$.
 \end{thm}
 
 {\it Proof:} Since $S$ is comprehensive the set $S\cap\mathcal{E}(u; X)$ is not empty, and it is a bounded chain in $S$. We have 
 $\sup[S\cap\mathcal{E}(u; X)]\in S$. Let $\bar{x} = \sup[S\cap\mathcal{E}(u; X)]$. We show that $u^\circ(\bar{x}) = \bar{x}$ since $\bar{x}$ is then the  largest efficient point of $S$ and by Proposition \ref{argmax1} 
 $\bar{x}\in\arg\!\max(u; S)$.\\ 
 For all $x\in S\cap\mathcal{E}(u; X)$ we have $x\leq  \bar{x}$  therefore  $u^\circ(x)\leq  u^\circ(\bar{x})$; since $u^\circ(x) = x$ if $x\in \mathcal{E}(u; X)$ we have shown that $u^\circ(\bar{x})$ is an upper bound of $S\cap\mathcal{E}(u; X)$ and therefore    $\bar{x}\leqslant  u^\circ(\bar{x})$. But we also have 
 $u^\circ(\bar{x})\geqslant \bar{x}$ and finally 
 $u^\circ(\bar{x}) = \bar{x}$. We have shown that $\arg\!\max(u; S)\neq\emptyset$. 
 
 \medskip\noindent Let $\bar{x}$ be an element of $\arg\!\max(u; S)$. If $\hat{x}$ if a maximal element of $S$ such $\bar{x}\leqslant\hat{x}$ then $u(\bar{x})\leqslant u(\hat{x})$ and therefore $\hat{x}\in \arg\!\max(u; S)$. We show that such an $\hat{x}$ exists.\\
 Let $x^\star$ be an upper bound of $S$ and consider the set $[\bar{x}, x^\star]\, \cap\, S$. Since $S$ has Property CUC we can invoke Zorn's Lemma to conclude that $\boldsymbol{Max}([\bar{x}, x^\star]\, \cap\,S)\neq\emptyset$. Obviously  
 $\boldsymbol{Max}([\bar{x}, x^\star]\, \cap\,S)\, \subset \, \boldsymbol{Max}(S)$.  \hfill$\Box$

 \begin{lem}\label{cuccomp} Assume that $X$ is a semilattice endowed with a topology for which intervals are compact and  the partial order is upper semicontinuous. If $S$ is comprehensive subset of $X$ such that, for all $x\in S$, $\uparrow\!\!(x)\cap S$ is closed  then $S$ is   (CUC).
 \end{lem} 
 
 {\it Proof:} Let $\mathcal{C}$ be a bounded chain in $S$ and let $\bar{x}\in X$ be an upper bound of $\mathcal{C}$. For all $x\in \mathcal{C}$, let 
 $U_x = [x, \bar{x}]\cap S$. Notice that $U_x$ is not empty since it contains $x$; furthermore, $U_x = [x, \bar{x}]\cap \big(\uparrow\!\!(x)\cap S \big)$, it therefore compact. For all finite subfamily $\{x_1, \cdots, x_m\}$ of the chain $\mathcal{C}$ there exists $i_0\in [m]$such that, for all $i\in [m]$, 
 $x_{i_0}\geqslant x_{i}$ and therefore $U_{x_{i_0}} = \cap_{i\in [m] }U_{x_{i}}$; this shows that $\{U_x: x\in \mathcal{C}\}$ is a family of compact sets with the finite intersection property. 
 Therefore $\cap_{ x\in \mathcal{C}}U_x\neq\emptyset$. An element of 
 $\cap_{ x\in \mathcal{C}}U_x$ is clearly an upper bound of $\mathcal{C}$ that belongs to $S$.\\
 Let $M$ be the set of upperbounds of $\mathcal{C}$ that belongs to $S$. We have seen that for all upper bounds $\bar{x}$ of $\mathcal{C}$ there exists $\tilde{x}\in M$ such that $\tilde{x}\leqslant \bar{x}$. We show that that $M$ has a smallest element. \\
 First, let us see that the set of minimal elements of $M$ is not empty. Let 
 $\mathcal{D}$ be a chain in $M$. For all $x_1, \cdots, x_m\in\mathcal{C}$ and all 
 $y_1, \cdots, y_n\in\mathcal{D}$, $\cap_{i\in[m], j\in[n]}[x_i, y_j] = [x_{i_0}, 
 y_{j_0}]$ where $x_{i_0} = \max\{x_1, \cdots, x_m\}$ and  
 $y_{j_0} = \max\{y_1, \cdots, y_n\}$; this shows that the family of compact sets $\{[x, y]: (x, y)\in\mathcal{C}\times\mathcal{D}\}$ has the finite intersection property. Let $x^\star\in\cap_{(x, y)\in\mathcal{C}\times\mathcal{D}}[x, y]$; $x^\star$ is a lower bound of $\mathcal{D}$ and an upper bound of $\mathcal{C}$. By Zorn's Lemma the set of minimal elements of $M$ is not empty. \\
 Let $M_0$ be the set of minimal elements of $M$. If $x_1, x_2\in M_0$ then, for all $x\in\mathcal{C}$, $x_1\geqslant x$ and 
 $x_2\geqslant x$ and therefore $ x_1\wedge x_2\geqslant x$. This shows that 
 $x_1\wedge x_2\in M_0$ and therefore $x_1 = x_1\wedge x_2 = x_2$. We have shown that $M_0$ contains  a single point $x^\star$. Obviously, $x^\star$ is the least upper bound of  $\mathcal{C}$ and $x^\star\in S$. \hfill$\Box$

\begin{prop}\label{cucuccomp} Let $X$ be a topological inf-semilattice  for which intervals are compact and let  $S$ be a   comprehensive and bounded above  subset of $X$ such that, for all $x\in S$, $\uparrow\!\!(x)\cap S$ is closed.  Then, for all  quasi-Leontief functions $u : X\to\Lambda$, $\arg\!\max(u; S)\neq\emptyset$. Furthermore,  
\begin{equation}\label{furthermore}\uparrow\!\!\big(\arg\!\max(u; S)\, \big) \cap S = \arg\!\max(u; S).
\end{equation}
\end{prop} 

{\it Proof:} From Theorem \ref{maxcuc} and Lemma \ref{cuccomp} we have $\arg\!\max(u; S)\neq\emptyset$; (\ref{furthermore}) holds because $u$ is isotone.\hfill$\Box$

\section{Individually quasi-Leontief functions}\label{secindiql}
In this section we consider a finite family of partially ordered set  
$\big(X_i, \leqslant_i\big)$, $i\in[n]$ and functions $u : \prod_{i\in[n]}X_i\to\Lambda$. \\ As usual, given $x= (x_1, \cdots, x_n)\in \prod_{i\in[n]}X_i$ we will denote by 
$x_{-j}$ the element of $\prod_{i\in[n]\!\setminus\{j\}}X_i$ obtained from $x$ by deleting $x_j$, we will also write $x = (x_{-j}, x_j)$;  we will also use the same notation for arbitrary elements of $\prod_{i\in[n]\!\setminus\{j\}}$ eventhough it might not  to be ``the $x_{-j}$ of a given $x$'' of the product. The partial order on $\prod_{i\in[n]\!\setminus\{j\}}X_i$ is the product partial order, that is 
$(x_1, \cdots, x_n)\leqslant (y_1, \cdots, y_n)$ if, for all $i\in [n]$, 
$x_i\leqslant y_i$. We will say that $u: \prod_{i\in[n]}X_i\to\Lambda$ is a {\bf globally quasi-Leontief function} if it is quasi-Leontief.\\
Given  $x_{-j}\in \prod_{i\in[n]\!\setminus\{j\}}X_i$ we denote by $u[x_{-j}]$ the function from $X_j$ to $\Lambda$ defined by 
$u[x_{-j}](x_j) = u(x_{-j}, x_j)$.

\begin{defn}\label{indiviquasileon} A function 
$u: \prod_{i\in[n]\!\setminus\{j\}}X_i\to\Lambda$ is individually (regular) quasi-Leontief if, for all $j\in[n]$ and all $x_{-j}\in \prod_{i\in[n]\!\setminus\{j\}}X_i$ the function $u[x_{-j}]$ is (regular) quasi-Leontief.
\end{defn}

\medskip\noindent We have seen in $(6)$ of Examples 
  \ref{exquaileon2} that globally (regular) quasi-Leontief functions are  
   individually (regular) quasi-Leontief.

\begin{expl}[Examples of individually quasi-Leontief functions]\label{exindivquaileon1}\quad\hfill\\
$\boldsymbol{[1]}$ $X_1 = X_2=\Lambda = \mathbb{R}_+$ and 
$u(x_1, x_2) = \min\{x_1, x_1x_2\}$. \\
$(a)$ $u$ is not globally quasi-Leontieff. Indeed, 
$u(x_1, x_2)\geqslant\lambda$ if and only if $x_1\geqslant\lambda$ and 
$x_1x_2\geqslant\lambda$ we can therefore take $x_1$ arbitrary large and 
$x_2$ small enough to have the second inequality; if $\lambda\neq 0$ the set $u^{-1}\big(\uparrow\!\!\lambda\big)$ does not have a smallest element. \\
$(b)$ $u$ is individually quasi-Leontief since, for all $(x_1, x_2)$, both $u[x_1](x) = \min\{x_1, x_1x\}$ and 
$u[x_2](x) = \min\{x, xx_2\}$ are minima of regular  quasi-Leontief functions on $\mathbb{R}_+$ ($4$ of Examples 
  \ref{exquaileon2}) . One can check that 
\\
$
u[x_1]^{-1}\big( u[x_1](x)\big) = \left\{
\begin{array}{lcl}
\, [x, \infty[ & \hbox{if }&x\leqslant 1    \\
\, [1, \infty[  &  \hbox{if}&x\geqslant 1    
  \end{array}
\right.
$

\medskip\noindent 
   $\boldsymbol{[2]}$ If $u_i: X_i\to\mathbb{R}_{+ +}$, $i\in [n]$, are quasi-Leontief functions then $u(x_1, \cdots, x_n) = \prod_{i\in[n]}u(x_i)$ is individually quasi-Leontief. (on could take $\mathbb{R}_{+}$ instead of $\mathbb{R}_{+ +}$ if one assumes that each $X_i$ has a smallest element to avoid potential difficulties with constant functions).
   
   \medskip\noindent 
   $\boldsymbol{[3]}$ If $u_i: X_i\to\mathbb{R}_{++}$, $i\in [n]$, are (regular) quasi-Leontief functions then, for all $(a_1, \cdots, a_n)\in\mathbb{R}_{++}$ and all $b\in\mathbb{R}$, $u(x_1, \cdots, x_n) = \sum_{i\in[n]}a_iu_i(x_i) + b$ is individually quasi-Leontief if each $X_i$ has a smallest element. 
  
     \medskip\noindent 
   $\boldsymbol{[4]}$ If $u_i: \prod_{j\in [m]}X_j\to\mathbb{R}$, $i\in [n]$, are regular individually quasi-Leontief functions and if each $X_j$ is a lattice then 
   $u = \min\{u_1, \cdots, u_n\}$ is regular individually quasi-Leontief.
 \end{expl}
 \subsection{Efficient points}\label{subseceff}
  
  The class of individually quasi-Leontief functions is so large that one cannot expect to have the kind of existence and uniqueness of efficient points that is characteristic of quasi-Leontief functions. For globally quasi-Leontief functions one can consider the set of efficient points or the set of points which are coordinatewise efficient given that the remaining coordinates are frozen; we will see that one can be recovered from the other.  
 
 \subsubsection{Efficient points for globally quasi-Leontief functions } 
  \begin{prop}\label{globeff}
  If $u : \prod_{i\in[n]}X_i\to\Lambda$ is globally quasi-Leontief then 
  $x\in\mathcal{E}\big(\prod_{i\in[n]}X_i, u\big)$ if and only if, for all 
  $j\in[n]$, $x_j\in\mathcal{E}(u[x_{-j}], X_j)$.\\
  Furthermore, for all $x\in\prod_{i\in[n]}X_i$ and for all $j\in [n]$, 
  $u^\circ(x)_j = u[x_{-j}]^\circ(x_j)$.
  \end{prop} 
  {\it Proof:} If $x\in \mathcal{E}\big(\prod_{i\in[n]}X_i, u\big)$ then 
  $u[x_{-j}](x_j^\prime) \geqslant u[x_{-j}](x_j) $ if and only if 
  $(x_{-j}; x_j^\prime)\geqslant (x_{-j}; x_j)$ if and only if  
  $ x_j^\prime\geqslant  x_j$ and therefore $x_j\in\mathcal{E}(u[x_{-j}], X_j)$. 
  Assume now that, for all $j\in [n]$, $x_j\in\mathcal{E}(u[x_{-j}], X_j)$. Since 
  $u$ is quasi-Leontief $u^\circ(x)$ is well defined and $x\geqslant u^\circ(x)$, in particular $x_j\geqslant u^\circ(x)_j $ and $x_{-j}\geqslant u^\circ(x)_{-j} $, and therefore 
  $x \geqslant (x_{-j}; u^\circ(x)_j)$ and also $(x_{-j}; u^\circ(x)_j)\geqslant u^\circ(x)$  from which we obtain 
  $u[x_{-j}](x_j)\geqslant u[x_{-j}](u^\circ(x)_j)\geqslant u\big(u^\circ(x)\big) = u(x)$. We have shown that $u[x_{-j}](x_j) =   u[x_{-j}](u^\circ(x)_j)$ and since 
   $x_j\in\mathcal{E}(u[x_{-j}], X_j)$ we must have $u^\circ(x)_j \geqslant x_j$. We have shown that  $x_j = u^\circ(x)_j $ or that $x = u^\circ(x)$. To prove the last part we have to show that $u[x_{-j}](x_j^\prime)\geqslant u[x_{-j}](x_j)$ if and only if $u[x_{-j}](x_j^\prime)\geqslant u^\circ(x)_j$ which is done exactly as the previous part. \hfill$\Box$
   
   \begin{prop}\label{globlambda} Assume that $u : \prod_{i\in[n]}X_i\to\Lambda$ is globally quasi-Leontief and regular. \\For all $\lambda\in\Lambda$ such that $u^{-1}\big(\uparrow\!\!(\lambda)\big)\neq\emptyset$ let 
   $u^\sharp(\lambda) = (u_1^\sharp(\lambda), \cdots, u_n^\sharp(\lambda))$. Then, for all $x \in \prod_{i\in[n]}X_i$ and for all $j\in [n]$, $u[x_{-j}]$ is regular quasi-Leontief and $u[x_{-j}]^\sharp(\lambda) = u_j^\sharp(\lambda)$  for all $\lambda$ such that $x_{-j}\geqslant u^\sharp(\lambda)_{-j}$  and therefore 
   $u_j^\sharp\big(u(x)\big) = u_j^\sharp\big(u[x_{-j}](x_j)\big) = u[x_{-j}]^\circ(x_j)$.
   \end{prop}
   {\it Proof:} From $u[x_{-j}](x_j)\geqslant\lambda$ if and only if 
   $(x_{-j}; x_j)\geqslant u^\sharp(\lambda)$ we have\\ 
   $\{x_j\in X_j: u[x_{-j}](x_j)\geqslant\lambda\}\neq\emptyset$ if and only if 
   $x_{-j}\geqslant u^\sharp(\lambda)_{-j}$ and $x_j\geqslant u^\sharp(\lambda)_{j}$; in other words, if $\{x_j\in X_j: u[x_{-j}](x_j)\geqslant\lambda\}\neq\emptyset$  then it has a smallest element, namely  
   $u_{j}^\sharp(\lambda)$. \hfill$\Box$
   
   \medskip\noindent
   From Proposition \ref{globlambda} it looks as if is  $u[x_{-j}]^\sharp$ does dot depend on $x_{-j}$; but one has to be carefull: 
   
   \medskip\noindent if $x_{-j}$, $x_{-j}^\prime$ and $\lambda$ are such that  $\{y\in X_j : u(x_{-j}; y)\geqslant\lambda\}\neq\emptyset$ and \\$\{z\in X_j : 
  u(x^\prime_{-j}; z)\geqslant\lambda\}\neq\emptyset$ then $u[x_{-j}]^\sharp(\lambda) = u[x^\prime_{-j}]^\sharp(\lambda) = p_j\circ u^\sharp(\lambda) = u^\sharp_j(\lambda)$ where $p_j$ is the projection of the product space $X$ onto $X_j$.  
   
   \medskip\noindent In particular, if $u[x_{-j}]^\sharp(\lambda)$ is defined and if $x_{-j}^\prime\geqslant x_{-j}$ then 
   $u[x^\prime_{-j}]^\sharp(\lambda)$ is defined and $u[x_{-j}]^\sharp(\lambda) = u[x^\prime_{-j}]^\sharp(\lambda)$; more generally, if 
   $u[x_{-j}]^\sharp(\lambda)$  and  
   $u[x^\prime_{-j}]^\sharp(\lambda)$ are both  defined then they are equal. For a simple example consider the function $u(x_1, x_2) = \min\{x_1, x_2\}$ defined on $\mathbb{R}_{++}^2$; then $u[a]^\sharp(\lambda) = \min\{x_2\in \mathbb{R}_{++}^2: u(a, x_2)\geqslant\lambda\}$ if defined only if $a\geqslant\lambda$ in which case it is $\lambda$.

   \bigskip\noindent
   \subsubsection{Efficient points for individually quasi-Leontief functions } 
    We know from Proposition \ref{globeff} that for a globally quasi-Leontief function the set $\mathcal{E}\big(\prod_{i\in[n]}X_i, u\big)$ of 
    efficient points is  the  fixed point set of the multivalued map defined on $\prod_{i\in[n]}X_i$ by 
    \begin{equation*}\mathbb{P}_u(x) = \prod_{i\in[n]}\mathcal{E}(u[x_{-i}], X_i).
    \end{equation*} Therefore, for a globally quasi-Leontief function, $x\in \mathbb{P}_u(x)$ if and only if $x$ is the smallest element, that is the unic minimal element, of $u^{-1}\big(\uparrow\!\!u(x) \big)$. 
       
 \noindent   The definition of $\mathbb{P}_u(x)$ makes sense for individually quasi-Leontief functions $u : \prod_{i\in[n]}X_i\to\Lambda$ but there is no reason for an arbitrary set $u^{-1}\big(\uparrow\!\!u(x) \big)$ to have a smallest element. The question then is what remains of Proposition \ref{globeff}  for individually quasi-Leontief functions. 
    
     \medskip\noindent
   First, the set  $\mathbb{P}_u(x)$ can be so large as to be of no interest. Indeed, let $X_1 = X_2 = \mathbb{R}_{++}$ and let $u(x_1, x_2) = x_1^{\alpha_1} x_2^{\alpha_2}$ be a Cobb-Douglas function with $\alpha_i> 0$; then 
   $u[x_2](x_1^\prime)\geqslant u[x_2](x_1)$ if and only if 
   $x_1^\prime\geqslant x_1$ and similarly $u[x_1](x_2^\prime)\geqslant u[x_1](x_2)$ if and only if 
   $x_2^\prime\geqslant x_2$; we have 
   $\mathcal{E}(u[x_{1}],  \mathbb{R}_{++}) =  \mathbb{R}_{++} = \mathcal{E}(u[x_{2}],  \mathbb{R}_{++})$ and therefore 
    $\mathbb{P}_u(x) = X_1\times X_2$. This happens because $X_i$ is one dimensional and $u[x_{-i}]$ is strictly increasing on $X_i$.

   \begin{defn}\label{pareff} Given a function  $u : \prod_{i\in[n]}X_i\to\Lambda$ we will say that a point $x\in\prod_{i\in[n]}X_i$ is  efficient if it is a minimal point of the set $u^{-1}\big(\uparrow\!\!u(x) \big)$ that is : 
   \begin{equation*}
   \forall x^\prime\,\,  \left[u(x^\prime)\geqslant u(x) \hbox{ and } 
   x\geqslant x^\prime\right] \Rightarrow x^\prime = x. 
   \end{equation*}
      \end{defn}
      
      \bigskip\noindent
      The set of efficient points can also be so large as being of no interest; with  the Cobb-Douglas function from the previous example one can easily see that $(x_1, x_2)\geqslant (x_1^\prime, x_2^\prime)$ and 
      $x_1^{\prime\alpha_1}x_2^{\prime\alpha_2}\geqslant x_1^{\alpha_1}x_2^{\alpha_2}$ implies 
      $(x_1, x_2) = (x_1^\prime, x_2^\prime)$. Therefore all points are efficient. 
      
       \medskip\noindent In this example we trivially have $x\in\mathbb{P}_u(x)$ if and only if $x$ is efficient, even if in this case it is a relatively uninteresting piece of information.  
       
   \medskip\noindent    We will see below that the equality between the fixed point set of  $\mathbb{P}_u$ and the set of efficient points always holds. But the situation does not have to be always as trivial as in the examples above. 
   
   \medskip\noindent For example, let 
       $X_1 = \mathbb{R}_{++}\times\mathbb{R}_{++}$ and $X_2 = \mathbb{R}_{++}$ with 
       $u(x_1, x_2, x_3) = \min\{x_1x_3, x_2\}$ which is individually quasi-Leontief since, for $(x_1, x_2) = (a, b)$ the partial function on $X_3$ is $u[(a, b)](x_3) =  \min\{ax_3, b\}$ which is quasi-Leontief on $\mathbb{R}_{++}$ and for $x_3 = c$ the partial function is $u[c](x_1, x_2) =  \min\{cx_1, x_2\}$ which is Leontief on 
       $ \mathbb{R}_{++}\times\mathbb{R}_{++}$. Let us see that $(x_1, x_2, x_3)$ is efficient if and only if 
       $x_1x_3 = x_2$. \\If $x_2 > \min\{x_1x_3, x_2\} = x_1x_3$ choose $(x_1^\prime, x_2^\prime, x_3^\prime)$ such that $x_1^\prime = x_1$, $x_3^\prime = x_3$ and $ x_2 > x_2^\prime > \min\{x_1x_3, x_2\}$; then 
       $(x_1, x_2, x_3)\geqslant (x_1^\prime, x_2^\prime, x_3^\prime)$,   
       $(x_1, x_2, x_3)\neq (x_1^\prime, x_2^\prime, x_3^\prime)$ and  
       $u(x_1, x_2, x_3) = u(x_1^\prime, x_2^\prime, x_3^\prime)$.\\
       If $x_1x_3 > x_2$ one proceeds similarly. If $x_1x_3 = x_2 =  \min\{x_1^\prime x_3^\prime, x_2^\prime\}$ with 
       $x_i\geqslant x_i^\prime$ then $x_i = x_i^\prime$.

   \bigskip
   
   \begin{prop}\label{charpar} Let  $u : \prod_{i\in[n]}X_i\to\Lambda$ be an individually quasi-Leontief function. For all $x\in \prod_{i\in[n]}X_i$ let 
   $\mathbb{P}_u(x) = \prod_{i\in[n]}\mathcal{E}(u[x_{-i}], X_i)$. Then a point 
   $x$ is  efficient if and only if $x\in \mathbb{P}_u(x)$.
   \end{prop}
   
    {\it Proof:} Assume that $x$ is  efficient. Let $x_j^\prime\in X_j$ be such that $u[x_{-j}](x_j^\prime)\geqslant u[x_{-j}](x_j)$. Since $u[x_{-j}]$ is quasi-Leontief there exists $x_j^{\prime\prime}\in X_j$ such that 
    $x_j\geqslant x_j^{\prime\prime}$, $x_j^\prime\geqslant x_j^{\prime\prime}$ and $u[x_{-j}](x_j^{\prime\prime}) = \min\{u[x_{-j}](x_j), u[x_{-j}](x_j^\prime)\}$, (Property $\Phi$), that is $u[x_{-j}](x_j^{\prime\prime}) = u[x_{-j}](x_j)$. We have $u(x_{-j}; x_j^{\prime\prime}) = u(x)$ and $(x_{-j}; x_j)\geqslant (x_{-j}; x_j^{\prime\prime})$; since $x$ is  efficient we must have 
    $(x_{-j}; x_j)= (x_{-j}; x_j^{\prime\prime})$, that is 
     $x_j = x_j^{\prime\prime}$ and therefore $x_j^\prime\geqslant x_j$. We have shown that $x_j\in \mathcal{E}(u[x_{-j}], X_j)$.
     
     \medskip\noindent
     Assume now that $x\in\mathbb{P}_u(x)$ and let $x^\prime$ be such that 
     $x\geqslant x^\prime$ and $u(x^\prime)\geqslant u(x)$. From $x\geqslant x^\prime$ we have $(x_{-j}; x^\prime_j)\geqslant x^\prime$ and therefore 
     $u(x_{-j}; x^\prime_j)\geqslant u(x^\prime)$ or, equivalently, 
     \begin{equation}\label{xprimej}
     x^\prime_j\geqslant u[x_{-j}]^\sharp\big(u(x^\prime)\big).
     \end{equation} 
     From $u(x^\prime)\geqslant u(x)$ we obtain 
     $u[x_{-j}]^\sharp\big(u(x^\prime)\big)\geqslant u[x_{-j}]^\sharp\big(u(x)\big)$ and therefore, from (\ref{xprimej}), 
     \begin{equation}\label{xprimej2}
     x^\prime_j\geqslant u[x_{-j}]^\sharp\big(u[x_{-j}])(x_j)\big).
     \end{equation}
     It follows from  $x_j\in \mathcal{E}(u[x_{-j}], X_j)$ that 
     $u[x_{-j}]^\sharp\big(u[x_{-j}])(x_j)\big) = x_j$ and consequently that 
     $x^\prime_j\geqslant x_j$. We have shown that $x^\prime\geqslant x$. \hfill$\Box$
     
     \begin{thm}\label{existseffglob} Let  $u : \prod_{i\in[n]}X_i\to\Lambda$ be an individually quasi-Leontief function and let $S_i\subset X_i$, $i\in [n]$ be comprehensive subsets.\\ If $\arg\!\max\big(u; \prod_{[i\in n]}S_i\big)\neq\emptyset$ then, for all $x^\star\in\arg\!\max\big(u; \prod_{[i\in n]}S_i\big)$ there exists 
     $x^{\flat}\in\arg\!\max\big(u; \prod_{[i\in n]}S_i\big)$ suchs $x^\flat$ is efficient and $x^\star\geqslant x^\flat$.
     \end{thm}
     
     {\it Proof:} Let $(x_1^\star, \cdots, x_n^\star) = x^\star\in\arg\!\max(u; S)$, where $S = \prod_{[i\in n]}S_i$, and let $x_1^\flat = u[x^\star_{-1}]^\circ(x_1^\star)$. From $u[x^\star_{-1}](x_1^\star) = u[x^\star_{-1}]\big(u[x^\star_{-1}]^\circ(x_1^\star)\big)$ we have 
     \begin{equation*}x^{[1]} = (x_1^\flat, x_2^\star \cdots, x_n^\star) \in\arg\!\max(u; S)
     \end{equation*}
     and from $x^\star_{-1} = x^{[1]}_{-1}$ we have 
     \begin{equation*}
     x_1^\flat\in \mathcal{E}(u[x^{[1]}_{-1}]; S).
     \end{equation*}
     Since $x_1^\star\geqslant x_1^\flat$ we also have $x^\star \geqslant x^{[1]}$.
     
   \medskip\noindent  Let $1\leqslant k\in[n-1]$ and assume that we have constructed a point \\ $x^{[k]} = (x_1^\flat, \cdots, x_k^\flat, x_{k+1}^\star, \cdots,  
   x_{n}^\star)\in S$ such that: 
   \begin{eqnarray}
 \label{ck1}  x^{[k]}\in \arg\!\max(u; S)\\\nonumber \\
 \label{ck2}  \forall i\in[k]\quad x_{i}^\flat\in \mathcal{E}(u[x^{[k]}_{-i}]; S_i)\\\nonumber \\
 \label{ck3}  x^\star\geqslant x^{[k]}.
   \end{eqnarray}
   
   \medskip\noindent Let $x_{k+1}^\flat = u[x^{[k]}_{-(k+1)}]^\circ(x_{k+1}^\star)$ and 
   $x^{[k+1]} = (x^{[k]}_{-(k+1)}; x_{k+1}^\flat)$. As in the first part of the proof one sees that $x^{[k+1]}\in \arg\!\max(u; S)$ and $x^{[k]}\geqslant x^{[k+1]}$, and therefore $x^{\star}\geqslant x^{[k+1]}$. From 
   $x^{[k]}_{-(k+1)} = x^{[k+1]}_{-(k+1)}$ we have 
   $x_{k+1}^\flat\in  \mathcal{E}(u[x^{[k+1]}_{-(k+1)}]; S_{k+1})$.\\
   Let us see that, for $i\in [k]$, $x_{i}^\flat\in  \mathcal{E}(u[x^{[k+1]}_{-i}]; S_i)$. We already have $x_{i}^\flat\in S_i$. Let $z_i\in X_i$ such that 
   $u[x^{[k+1]}_{-i}](z_i)\geqslant u[x^{[k+1]}_{-i}](x^\flat_i) = u(x^{[k+1]}) $. Since, for all $i\in[k+1]$, $x^{[i]}\in \arg\!\max(u; S)$, we can write  
   \begin{equation}\label{uk+1}
   u[x^{[k+1]}_{-i}](z_i)\geqslant  u(x^{[i]}).
   \end{equation} 
   From $(x_{i+1}^\star, \cdots, x_n^\star)\geqslant (x_{i+1}^\flat, \cdots, 
   x_{k+1}^\flat, x_{k+2}^\star, \cdots,  x_n^\star) $ we have 
   \begin{equation}\label{uk+2} 
   u[x^{[i]}_{-i}](z_i) = u(x_1^\flat, \cdots, x_{i-1}^\flat, z_i, x_{i+1}^\star, \cdots, x_n^\star)\geqslant u[x^{[k+1]}_{-i}](z_i).
   \end{equation}
   From (\ref{uk+1}) and (\ref{uk+2}) we obtain $u[x^{[i]}_{-i}](z_i)\geqslant  u(x^{[i]})$ and, since $x_{i}^\flat\in \mathcal{E}(u[x^{[i]}_{-i}]; S_i)$, 
   $z_i\geqslant x_i^\flat$.  We have shown that (\ref{ck1}), (\ref{ck2}), and 
   (\ref{ck3}), hold for $x^{[k+1]}$. In a finite number of steps we obtain a point $x^\flat = (x_1^\flat, \cdots, x_n^\flat)$ such that $x^\flat\in\mathbb{P}_u(x^\flat)$, $x^\flat\in\arg\!\max(u, S)$ and $x^\star \geqslant x^\flat$. By Proposition \ref{charpar}  $x^\flat$ is efficient.  
   \hfill$\Box$


\begin{thebibliography}{99}
\bibitem{alibor}{\sc Aliprantis C. D. and Border K.}, {\it Infinite Dimensional Analysis: A Hitchhiker's Guide },  SpringerVerlag, 3rd. edition, 2006.

\bibitem{agk}{\sc Allamigeon X., Gaubert S., Katz D.}{\it Tropical Polar Cones, Hypergraph transversals, and mean payoff games}, arXiv:1004.2778v2

\bibitem{blyth}{\sc Blyth, T.S.}, {\it Lattices and Ordered Algebraic Structures}, Springer, 2005.

\bibitem{blythjano}{\sc Blyth, T.S and Janowitz, M.F.}, {\it Residuation Theory}, Pergamon Press, 1972.


\bibitem{t98}{\sc Topkis, D.}, {\em  Supermodularity and Complementarity}, Princeton University Press, 1998.

\bibitem{t89}{\sc\sc Topkis, D.}, {\it Equilibrium points in nonzero-sum n-person submodular games}, SIAM
Journal of Control and Optimization, vol. 17, pp. 773-787.


\end{thebibliography}
\end{document}